\documentclass[a4paper,11pt]{article}
\usepackage[T1]{fontenc}
\usepackage[english]{babel}
\usepackage[latin1]{inputenc}
\usepackage{amsmath}
\usepackage{amsmath, amsthm}
\usepackage{amsfonts,amssymb}
\usepackage{mathrsfs}
\usepackage{cases}
\usepackage{color}
\usepackage{natbib}
\usepackage{comment}

\usepackage{empheq} 
\usepackage[top=3.5cm, bottom=4cm, left=2cm , right=2.1cm]{geometry}

\theoremstyle{plain}
\newtheorem{theo}{Theorem}[section]
\newtheorem{req}{Remark}[section]
\newtheorem{Req}{Remark}[section]
\newtheorem{pp}[theo]{Proposition}

\newtheorem{df}{Definition}[section]

\newtheorem{cor}[theo]{Corollary}

\numberwithin{equation}{section}

\date{November 30, 2012 }

\title{An Exact Connection between two Solvable SDEs and a  Nonlinear Utility Stochastic PDE  \footnote{With the financial
support of
the "Chaire Risque Financier of the  Fondation du Risque", the " Chaire D\'eriv\'es du futur" 
of the "F\'ed\'eration des banques Fran\c caises".} }

\author{ El Karoui Nicole,
\thanks{ \small  LPMA, UMR CNRS  6632,  Universit\'e Pierre et Marie Curie, CMAP, UMR CNRS 7641, \'Ecole Polytechnique  }
\\
\and Mrad~Mohamed~ \thanks
{\small CMAP, UMR CNRS 7641, \'Ecole Polytechnique, \small  LAGA, UMR CNRS 7539,  Universit\'e Paris 13}}

\def\ind{{\bf 1}}
\def\demi{\frac{1}{2}}

\def\cal{\mathcal}

\def\L{{\cal L}}

\def\F{{\mathcal F}}
\def\K{{\cal K}}

\def\Kc{{\cal K}}

\def\sigR{{\cal R}}
\def\sigoR{{\cal R}^{\Sigma,\perp}}
\def\sigoR{{\cal R}^{\perp}}
\def\sR{\text{{\tiny${\cal R}$}}}
\def\S{{\cal S}}
\def\X{{\cal X}}

\def\Y{{\cal Y}}
\def\Z{{\cal Z}}
\def\bfB{{\bf B}}
\def\bfF{{\bf F}}
\def\bfM{{\bf M}}

\def\bfU{{\bf U}}

\def\E{{\mathbb E}}
\def\bF{{\mathbb F}}

\def\P{{\mathbb P}}

\def\R{{\mathbb R}}

\def\R{{\mathbb R}}

\def\PX{{\mathbb X}^+ }

\def\oK{{\overline{\mathcal K}}}

\def\tU{{\widetilde U}}
\def\tu{{\tilde u}}
\def\tv{{\tilde v}}
\def\wL{{\widehat L}}
\def\wL{{\widehat L}}

\def\tV{{\widetilde V}}

\def\Cc{{\mathcal C}}

\def\GX{{\mathscr X}}
\def\GY{{\mathscr Y}}

\def\eps{\varepsilon}

\newcommand{\rmi}{{\rm (i) $\>\>$}}

\newcommand{\rmii}{{\rm (ii) $\hspace{1.5mm}$}}
\newcommand{\rmiii}{{\rm (iii)$\>\>$}}
\newcommand{\rmiv}{{\rm (iv)$\>\>$}}
\newcommand{\rmv}{{\rm (v)$\>\>$}}

\newcommand{\rma}{{\rm a)$\>\>$}}
\newcommand{\rmb}{{\rm b)$\>\>$}}
\newcommand{\rmc}{{\rm c)$\>\>$}}

\def\edoc { \end{document}}
\def\bcom{}

\begin{document}
 \maketitle
 \abstract{Motivated by the work of Musiela and 
Zariphopoulou \cite{zar-03}, we study the It\^o random fields which are utility functions  $U(t,x)$ for any $(\omega,t)$.
The main tool is the marginal utility $U_x(t,x)$  and its inverse expressed as the opposite of the derivative of the Fenchel
conjuguate  $\tU(t,y)$. 
Under regularity assumptions, we associate 
 a $SDE(\mu, \sigma)$ and its adjoint SPDE$(\mu, \sigma)$ in divergence form whose $U_x(t,x)$ and its inverse $-\tU_y(t,y)$ are
monotonic solutions. More generally, special attention is paid to rigorous justification of  the dynamics of inverse flow of SDE.
So that, we are able to extend  to the solution of  similar SPDEs  the decomposition based on  the solutions of two
SDEs and their inverses.

The second part is concerned with forward utilities, consistent with a given incomplete financial market, that can be observed but given 
exogenously to the investor. As in \cite{zar-03},  market dynamics are considered in an equilibrium state, so that the
investor  becomes indifferent to any action she can take in such a market.
After having made explicit the constraints induced  on the local characteristics of consistent utility and its
conjugate, we focus on the marginal utility SPDE  by showing that it belongs to the previous family of SPDEs.
The associated two SDE's are related to 
the optimal wealth and the optimal state price density, given a pathwise explicit representation of the marginal utility.
This new approach addresses several issues with a new perspective: dynamic programming  principle, risk
tolerance properties, inverse problems. Some examples and applications are given in the last section}

\section*{Introduction}
The concept
 of forward dynamic utility has been introduced by M. Musiela and T. Zariphopoulou (2003-2008 \cite{zar-03,zar-07a,
zar-08,zar-07}), to model possible changes over the time of individual preferences of an agent. 
 The economic agent will adjust its preferences  in uncertain universe based on the information that is revealed over time 
 and  represented by a filtration $\bF=(\F_t)_{t \ge 0}$, defined on the probability space $(\Omega,\P, \cal A)$. In contrast to the classical literature, there is no pre-specified trading horizon at the end of which the utility datum is assigned. 
 Rather, the agent
starts with today's specification of its utility, $U(0,x)=u(x)$, and then builds the
process $U(t,x)$ for $t > 0$ in relation to the information  given by $(\F_t)_{t \ge 0}$.
This, together with the choice of an initial utility, distinguishes  forward dynamic utility from
 recursive utility first introduced by Duffie and Epstein, \cite{Duffie}, (see also \cite{NEKPENG}) for which the aggregator can be specified exogenously and
the value function is recovered backward in time. 
Dynamic utility may be also defined in relation with a exogenous investment universe $\X$, described as an incomplete market whose the parameters are known by the investor. Following \cite{zar-03},  market dynamics are considered in an equilibrium state, so that the investor  becomes indifferent to any action she can take in such a market. She is forced to embed into his dynamic
utility, structure information describing the equilibrium market dynamics. Such $\X$- consistent utilities will be studied in
the second part of the paper.\\

Such  concept has also been
studied by Berrier, Tehranchi and Rogers (2009) \cite{Mike} and Zitkovic \cite{zitkovic}. Further works  related
to this problem are Choulli, Stricker and Jia (2007) \cite{choulli},   Henderson and Hobson (2007)
\cite{Hobson}.

We start  by introducing the general notion of progressive utility and its Fenchel transform.
However, we restrict our study to positive
progressive utility which is a collection of It\^o's semimartingales with respect to a $d$-dimensional Brownian motion $W$, depending on a spatial parameter $x \in \R^+$, also called the wealth. As usual, they are specified
through their local (drift, volatility) characteristics $(\beta, \gamma)$ and their initial condition, a deterministic concave function $u$,  by
\begin{equation}\label{dyn0}
 dU(t,x)=\beta(t,x)dt + \gamma(t,x).dW_t, \quad U(0,x)=u(x).
\end{equation}
Apart regularity issues, the first question to be solved is to
give conditions on the local characteristics $(\beta, \gamma)$, insuring the
properties of monotonicity and concavity of the progressive utility $\bfU$
defined by the forward equation \eqref{dynbetagamma}.
The problem is equivalent to show that the progressive marginal utility $\bfU_x$
is strictly decreasing and strictly positive, with range $(0,\infty)$. At this stage, it is convenient to introduce the Fenchel conjugate random field $\bf \tU$ of $\bfU$ whose the opposite of the derivative $\bf \tU_y $ is the inverse flow of $\bfU_x$. 
Following Kunita \cite{Kunita:01}, we introduce sufficient conditions implying that $\bfU_x$ is also an It\^o progressive random field, with local characteristics
$(\beta_x, \gamma_x)$. (Theorem \ref{DRules}. All regularity issues are presented in details in Section \ref{subsec:regularityKunita}.)

 Then we show that the random field $\bfU_x$ is a monotonic with respect to its initial condition, strong solution of a one dimensional
stochastic differential equation. The random coefficients of the SDE$(\mu,\sigma)$ are expressed in terms of the local
characteristics $ \beta_x $ and $\gamma_x$, and $-\tU_y$. We give in particular, under some Lipschitz
assumptions on   SDE's coefficients, a  clear answer to the existence,
monotonicity, regularity and concavity  of
random fields defined from \eqref{dyn0} satisfying the Inada's conditions,
which until now remained an open question.\\
Additional regularity are need to show that the conjugate random field  $\tU$ and its derivative $\tU_y  $ are also  families of It\^o's semimartingales, and to identify its local characteristics. The It\^o-Ventzel formula enables us to carry out computations in a stochastically modulated dynamic framework.
Applied to the random field $\bfU$ and the semimartingale $\tU_y(t,y)$, It\^o-Ventzel's formula provides us two complementary
results: the dynamics of $\tU$, and the characterization of $\tU_y(t,y) $ as monotonic solution of the SPDE generated by the
adjoint elliptic operator in divergence form associated with the SDE$(\mu,\sigma)$ satisfied by $\bfU_x$. A first connection
between SDE and SPDE is then established.\\
Section \ref{subsec:regularityKunita} is more technical. Based on the reference book  of Kunita \cite{Kunita:01}, we provide
additional results in both regularity of random fields, and behavior of SDE with local Lipschitz property. The aim is to
provide conditions under which the inverse of a monotonic solution of SDE$(\mu,\sigma)$ (for instance the process
$-\tU_y(t,y)) $ is an It\^o semimartingale, solution of the second-order SPDE$(\mu,\sigma)$ (Theorem \ref{SDEInv}).
Conversely,  a regular  monotonic solution of SPDE$(\mu,\sigma)$ is identified with the inverse of the monotonic solution of 
SDE$(\mu,\sigma)$.  The link between SDE and second order SPDE is extended to more complex SPDE in the main result  of this
section, Theorem \ref{pp:  composition formula}.

In the second part of the paper, we turn to consistent stochastic utilities, where  the preferences are no more defined in isolation with a financial
market. The dynamics of the market
can be observed but are exogenously given to the investor, modulo estimation of
the relevant parameters by the equilibrium considerations or other procedures.
Following \cite{zar-03},  market dynamics are considered in an equilibrium state, so that the investor  becomes indifferent
to any action she can take in such a market. But, once characterized, the criterion may be used to solve classical
optimization problems as optimal allocation, indifference pricing \cite{Rouge}.\\
Working with positive wealth processes $X^\kappa$ in an  incomplete market,  we
define a consistent stochastic utility as a progressive  non-negative stochastic utility $U(t,x)$, for which 
$U(t,X^\kappa_t)$ is a supermartingale, and a (local) martingale for a certain portfolio called optimal wealth, or optimal benchmark.
Section \ref{CDU} is dedicated to  set-up  the financial market, and to identify the constraints induced on the local
characteristics $(\beta, \gamma)$ of $\bfU$ by the consistency property, using once again the It\^o-Ventzel formula.
As in
the classical Hamilton-Jacobi-Bellman framework, we proceed by verification to 
establish the dynamics of consistent utilities. Assuming a sufficient 
constraint on the drift  $\beta$ of HJB type, we get the utility SPDE
that we investigate in this paper. In particular, we study the role of  the
utility risk premium  defined by  $\eta^{U}(t,x)=\gamma_x(t,x)/U_x(t,x)$.
Simultaneously, we go into details on the duality questions  and give a characterization of the 
nonlinear SPDE satisfied by the progressive convex conjugate $\tU$ of $U$, and its optimal coefficients. \\
Our most original contribution is provided in Section \ref{USPDEOSDE}, using a characterization of $\X$-consistent dynamic
utility from its marginal utility. Based on  ideas developed in Section \ref{subsec:regularityKunita}, we show that, under
regularity assumptions on the random field $\bfU$, $U_x$  is solution of  a SPDE\eqref{eq: composition formulaA} associated
with the optimal coefficients for the wealth and state price density processes. Then, a first result shows that under natural
assumptions yielding to the existence of strong optimal solution of the conjugate problem $Y^*_t(y)$, the only locally
Lipschitz SDE$(\mu^*,\sigma^*) $ has a monotonic non-explosive solution $X^*_t(x)$. A simple application of Theorem \ref{pp: 
composition formula} yields to the closed form of the marginal utility as $U_x(t,x)=Y^*_t(u_x((X^*_t)^{-1}(x)))$.\\
Then, we consider the converse problem of recovering a dynamic utility coherent with a given optimal portfolio. In the classical case, this problem is
known in the financial literature as
the ``inverse`` Merton problem; it has been considered by many authors  in particular  by H.He and C.Huang (1992),
\cite{HeHuang}. Since the class of dynamic utilities is larger than the class of Markovian utilities considered in  \cite{HeHuang}, our
problem is easier to solve. Given a monotonic admissible portfolio $\bf X$, for any {\bf regular} state price density process $Y_t(y)$ and 
any initial condition $u_x$ such that $V(t,x)=Y^*_t(u_x((X_t)^{-1}(x)))$ is integrable near to $0$, is the marginal utility of some consistent dynamic utility.
So, we are able to generate all the consistent utilities with a given optimal portfolio.\\
\noindent
In Section \ref{USPDEOSDE}, Theorem \ref{ppU_x}, we establish the most original contribution of this paper. Indeed, from
Theorem \ref{pp: composition formula}, focusing on the dynamics of the marginal utility $U_x$, we can read directly and
without difficulties that the process $U_x\big(t, X^*_t(x)\big) $ is equal to the optimal dual state price density process
$Y^*(t,u_x(x))$ with  $u$ is the initial utility function. From this dual identity, the idea is very simple and natural:
Suppose that the optimal portfolio denoted by $\big(X^*_t(x)\big)_t$ is  strictly increasing
 with respect to the initial capital, and denote by $\big(\X(t,x)\big)_t$ the adapted inverse process, defined by
$X^*_t\big(\X(t,x)\big)=x$ then, we can find $U_{x}(t,x)$ from
$U_{x}(t,x)=Y^*_t\big(u_x(\X(t,x))\big)$.
Finally we get $U$ by integration. So, we are able to generate all the consistent utilities with a given optimal portfolio.\\
We recover easily  the  martingale property of the risk tolerance established in 
\cite{HeHuang}, in a complete market.
 Finally, in Section  \ref{EDSU}, we close the paper by some openness to other topics and works;  we show  the stability of
the notion  of consistent utility by change of numeraire and then,  without loss of generality,  we can  consider the  martingale market where the
portfolios are simple local martingales and the stochastic PDE's are easier to deal with.  
We also apply our method to the specific
example of decreasing consistent utilities (see \cite{Mike}  and \cite
{zar-08a}) where the volatility vector $\gamma$ is given equal to zero, given a new interpretation of the optimal wealth as solution of sup-convolution problems in random power utilities.\\

To the best of our knowledge, the characterization of strictly concave
increasing random fields, satisfying Inada's conditions, from SDEs have neither
studied in details  nor established in the literature. Moreover, 
the  fully  nonlinear utility stochastic PDE's established in this paper and  satisfied by forward utilities and their dual
have not been   established in a general way. 
Furthermore,  there is no general consistent utilities construction proposed in the literature, expect the
case of power or exponential type, or decreasing utilities.
Another main contribution of this paper is  a connection between two solvable
SDEs and the utility SPDE early established. In particular, given a  volatility
vector $\gamma$ such that
$\gamma_x(t,x)=-xU_{xx}\kappa^*_t(x)-U_x(t,x)\eta^{\cal
R}_t+U_x(t,x)\nu^*_t(U_x(t ,x))$ (where $\eta^\sR$ is the market
risk prime), we show the existence and uniqueness of a solution to  the fully 
nonlinear second order SPDE from that of a pair of SDE's. In any case this,
represents an interesting result in the theory of stochastic partial
differential equations.\\

  The paper is organized as follows, in  the next section we give the definition
of the It\^o progressive  dynamic utilities and focus on the characterization of these concave
It\^o's random fields by establishing a link with SDEs. A special attention is paid to the dynamics of the Fenchel conjuguate utility random field. 
In the following section, we give sufficient conditions ensuring concavity, monotonicity, differentiability both for random fields or for solution of SDE.
The dynamics of the inverse of regular solution of SDE$(\mu, \sigma)$ is the first step in the links of SDE and SPDE. More complex situations are then studied.

 Next, we
introduce the consistent utilities, the market model and the framework of the
paper.  In Section  \ref{CDU}, we  provide, simultaneously, the dynamics of
consistent utilities and their convex conjugate. We establish a closed form for
the optimal policies.  Thereafter, in Section  \ref{USPDEOSDE}, we establish a fully characterization of the marginal
utility. We 
end by reverse engineering problem and some openness to other topics.

\section{Progressive Utility}\label{SPU}
In a dynamic and stochastic environment, the classical notion of utility is not
flexible enough  to help us to make good choices in the long run. M.Musiela and
T.Zariphopoulou  were
the first to suggest to use instead of the classical criterion  the concept of
 progressive dynamic utility, consistent with respect to a given investment universe in a sense specified in
Section
3.
The concept of progressive utility  gives an adaptative way to include  new
information on environment evolution available to economic agents. 
Recently, Fritelli \cite{Fritelli} introduced very closed
notion called stochastic dynamic utility, in view of study certainty equivalent. Since these
utility functions are stochastic, time dependent and moving forward, we
consider them as a family of (It\^o) semimartingales depending on a parameter, the
wealth of the agent in the economic context.
\subsection{Definition and Properties of Progressive Dynamic Utility}
All stochastic processes are defined on a standard filtered probability space
$(\Omega,{\mathbb F},\mathbb{P})$, where
the filtration ${\bF}=(\F_t)_{t\ge 0}$ is assumed to be right continuous and
complete.  We first recall some notions relative to stochastic processes
depending of a parameter.
\paragraph{Generality on progressive random fields}\label{par:PRFprop}
 Stochastic processes in consideration are   depending  on a parameter, often
called spatial parameter. For us, because of economic motivation, this parameter
is  the wealth of an investor, taking non negative values in $\R^*_+=\{x>0\}$. 
 Sometimes, we will use the vocabulary of random field theory, and refer to such
processes as progressive random fields. As all random fields considered in
the sequel are progressive, we  will often omit the mention " progressive".\\
\rmi A progressive random field ${\bf X}=\{X (t, x); t \geq 0, x > 0\}$ is  a
 random variable measurable with respect to $\F_\infty\otimes {\mathcal
B}(\R_+)\otimes {\mathcal B}(\R^*_+)$, which is a collection of progressive
processes $t\mapsto X(t, x)$.\\
\rmii A random field ${\bf X}$ is said to satisfy a property  $\cal P$ with
respect to $x$, if there exists $N\in \F_{\infty}$ with $\P(N)=0$, such that for
any $\omega\in N^c$, and any $t\geq 0$, $x\mapsto X(t,x)(\omega)$ satisfies the
property $\cal P$. For  example,  ${\bf X}$ is said to be concave, (resp.
increasing) if there exists $N\in \F_{\infty}$ with $\P(N)=0$, such that for any
$\omega\in N^c$, and any $t\geq 0$ $x\mapsto X(t,x)(\omega)$ is concave (resp.
increasing). \\
\rmiii A random field ${\bf X}$ is said to be {\em continuous}  if
there exists $N\in \F_{\infty}$ with $\P(N)=0$, such that for any $\omega\in
N^c$, and any $t\geq 0$ $x\mapsto X(t,x)(\omega)$ is continuous. In general, the existence of continuous modification is
based of the Kolmogorov's criterion, that we recall in Section \ref{subsec:regularityKunita}.\\
A random field ${\bf X}$ is said to be  differentiable if there exists $N\in \F_{\infty}$ with $\P(N)=0$, such that for any
$\omega\in
N^c$, and any $t\geq 0$ $x\mapsto X(t,x)(\omega)$ is differentiable
; the derivative denoted $X_x(t,x)(\omega)$ generates the so-called derivative random field ${\bf X}_x$.
When  ${\bf X}_x$  has a  continuous version, ${\bf X}$ is said to be ${\mathcal C}^1$-{\em regular}.\\
 \rmiv  A  d-dimensional random field ${\bf X}$ is said to be $\L^k$- {\em locally bounded}
 if for any compact $K\subset]0,\infty[$, there exists $N\in
\F_{\infty}$ with $\P(N)=0$, such that for any $\omega\in N^c$,  $\int_0^T\sup_{x\in K}\|X(t,x)\|^k(\omega)dt<\infty$  for
any $T>0$.

In the sequel, we are concerned with differentiable random fields said to be ${ \K}^m_{loc}$-regular ($\K$ for Kunita) in the
following sense.
\begin{df} \label{Kregularity} A  $d$-dimensional random field ${\bf X}$ is said to be ${\K}^m_{loc}$(resp.
${\oK}^m_{loc}$)-regular if $X$ is a $\Cc^m
$-regular random field such that  $X/x$, $\partial_x^{k}X,  k\leq m $  are $\L^1$(resp. $\L^2$)-locally bounded. Such random
fields are also called of class ${\K}^m_{loc}$(resp. ${ \oK}^m_{loc})$.

\end{df}
At this stage, we content ourselves with this definition to carry out calculations in this section; more details are given in
the next section.

\paragraph{Progressive utility and its Fenchel conjugate}
We start with the definition of a progressive utility as progressive random
field with  concavity property.
\begin{df}[{Progressive  Utility}]\label{defPUF}
A {\em progressive utility} is a continuous progressive random field on $\R^*_+$,
$\bfU=\{U(t, x); t \geq 0, x > 0\}$ such that,

\rmi {\sc Utility property:} $\bfU$ is strictly concave, strictly increasing, and 
non negative.
 
\rmii  {\sc Regularity property:} $\bf U$ is  a ${\mathcal C}^2$-random field,
with continuous 
 first and second derivatives random fields ${\bf U_x}$ and ${\bf U_{xx}}$.

\rmiii {\sc Inada conditions:} $\bf U$ goes to $0$ when $x$ goes to $0$ and the
derivative $\bf U_x$ goes to $\infty$ when $x$ goes to $0$, and to $0$ when $x$
goes to $\infty$.
\end{df}
Given its importance in convex analysis, we introduce together with any
progressive utility $\bf U$, its convex conjugate $\bf \widetilde U$ (also
called conjugate progressive  utility (CPU)), that is the Fenchel Legendre
transform of the convex random field $-\bfU(,-.)$. 

\begin{df}[{Progressive conjugate utility}]\label{defCPUF}
 The convex conjugate of the progressive utility $\bf U$ is
the progressive random field $\bf \tU$ defined on  $\R^*_+$ by 
${\bf \tU}=\{\tU(t,y); t\geq 0,y> 0\}$, where
$\tU(t,y)\stackrel{def}{=}\max_{x>0,x\in Q^+}\big(U(t,x)-x\, y\big)$.

\rmi Under  Inada condition, ${\bf \tU}$ is twice continuously
differentiable, strictly convex, strictly decreasing, with 
$\tU(.,0^+)=U(+\infty), \>
\tU(.,+\infty)=U(0^+), a.s.$

\rmii The marginal  utility $\bfU_x$  is the inverse of the opposite
of the   marginal conjugate utility ${\bf \tU}_y$, that is
$U_x(t,.)^{-1}(y)=-\tU_y(t,y)$, with Inada conditions $\tU_y(.,0^+)=-\infty$, 
$\tU_y(.,+\infty)=0$

\rmiii The bi-dual relation holds true $U(t,x)=\inf_{y>0,y\in
Q^+}\big(\tU(t,y)+x\, y\big)$. \\Moreover $\tU(t,y)=U\big(t,-\tU(t,y)\big) +\tU_y(t,y)\, y$,
and $U(t,x)=\tU\big(t, U_x(t,x)\big)+x\, U_x(t,x)$. 
\end{df}

\subsection{It\^o Progressive Utility and SDE}\label{IPU}
In this paper, we focus  on  continuous progressive utilities $\bfU$ which are 
a collection of  It\^o  semimartingales: for any $x$, $U(.,x)$ is a continuous
It\^o semimartingale, driven by a  $n$-dimensional Brownian
motion $W=(W^1,..,W^n)$ defined on the probability space $(\Omega,{\mathbb
F},\mathbb{P})$. In general, the $\sigma$-field $\F_0$ is assumed to be trivial
and $\F_0$- random variables are a.s. constant. We refer to the book of H.Kunita \cite{Kunita:01}  and to the next section
for all technical results
concerning the theory of semimartingale random fields.The assumption of finite
dimensional Brownian motion greatly simplifies  the theory.
 
 As usual, an It\^o random field $\bfF$ defined on $\R^+$ is specified through its decomposition $\bfF=F_0+\bfB^{F}+\bfM^F$
into two random fields, where $\bfB^{F}$ is a finite variation random field $B^F(t,x)=\int_0^t \varphi^F(s,x)ds$ and $\bfM^F$
is a martingale random field $M^F(t,x)=\int_0^t\psi^F(s,x)dW_s$, where $(\varphi^F, \psi^F)=\{(\varphi^F(t,x), \psi^F(t,x));
t\geq 0,x>0\}$ are
the local characteristics of $F$ assumed
to be progressive random fields, with values in $\R$ and  $\R^n$ respectively.  $ \varphi^F$ is called the drift
characteristic, and  $ \psi^F$ the
diffusion characteristic. For simplicity, we often omit the index $F$. 
By convention, an It\^o random field $\bfF$ is said to be a ${ \K}^m_{loc}$-{\bf semimartingale}, whenever $ \bfB^F$ is of
class ${ \K}^m_{loc}$, and $ \bfM^F$ is of class ${ \oK}^m_{loc}$. The reference to ${ \K}^m$ recall that $\bfF$ is a random
field.\\

Let us consider an It\^o progressive utility $U$ with initial condition $u(x)=U(0,x)$, and local characteristics $( \beta,
\gamma)$,
 so that
 \begin{equation}\label{dynbetagamma}
 dU(t,x)=\beta(t,x)dt + \gamma(t,x).dW_t,\quad U(0,x)=u(x)
\end{equation}
 Apart regularity issues, the first question before going further is the following:\\
{\em How to express by conditions on the local characteristics $(\beta, \gamma)$ random fields, the
properties of monotonicity and concavity of the progressive utility $\bfU$
defined by \eqref{dynbetagamma}?}\\
The problem is equivalent to show that the progressive marginal utility $\bfU_x$
is strictly decreasing and strictly positive, with range $(0,\infty)$.
Following Kunita \cite{Kunita:01},
this can be done  by assuming  that, \\
\rmi $U$ is of class ${\K}^2_{loc}$, which implies in particular (Theorem \ref{DRules} below) that $\bfU_x$
is
also an It\^o progressive random field, with local characteristics
$(\beta_x, \gamma_x)$. \\
\rmii The random field $\bfU_x$ is a strong solution of a one dimensional
stochastic differential equation (SDE), with random coefficient, monotonic with respect to its initial condition. \\
We summarize these ideas in the following theorem.

\begin{theo}\label{UxSDE} 
 Let  $\bfU$ be an  It\^o random field with dynamics 
 \begin{equation}
 d\,U(t,x)=\beta(t,x)dt + \gamma(t,x).dW_t,\quad U(0,x)=u(x)
\end{equation} 
\rmi {\sc Necessary condition} Let $\bfU$ be a progressive utility with conjugate
utility $\bf \tU$, and marginal utilities $U_x(t,.)$ and $-\tU_y(t,.)$.
If $\bfU$ is a  ${\K}^2_{loc}$-semimartingale, the random field $Z_.(z)=U_x(.,-\tu_y(z))$
is a strictly increasing (in $z$) solution of the SDE$(\mu, \sigma)$, 
\begin{equation}\label{eq:ItoSDE}
\left\{
\begin{array}{cllll}
 dZ_t&=\mu(t,Z_t)dt+\sigma(t,Z_t)\,dW_t, & \quad Z_0 = z\\
\mu(t,z)&=\beta_x\big(t,-\tU_y(t,z)\big),& \quad\mu(t,0)=0\\
\sigma(t,z)&:=\gamma_x\big(t,-\tU_y(t,z)\big),&\quad \sigma(t,0)=0
  \end{array}
 \right .
\end{equation}

\noindent
\rmii {\sc Characterization as primitive of SDE}  Let consider  
a SDE$(\hat{\mu},\hat{\sigma})$, $dZ_t=\hat{\mu}(t,Z_t)dt+\hat{\sigma}(t,Z_t)\,dW_t,  \> Z_0 = z$ and assume the existence of
 a strong solution $Z_.(z)$,  increasing and differentiable
 in $z$ with range $(0,\infty)$. 
Then, for any utility function $u$ such that $Z_.(u_x(x))$ is integrable in a
neighborhood of $x=0$, 
the primitive $\bfU=\{U(t,x)=\int_0^xZ_t(u_x(z))dz, t \geq 0, x>0\}$ is a
progressive utility.
\end{theo}
\begin{proof}
 By assumption, from  Theorem \ref{DRules} below, $U$ is a ${ \K}^2_{loc}$-semimartingale implies that $U_x$ is an It\^o
semimartingale with local characteristics $(\beta_x,\gamma_x)$, i.e.,   
$ dU_x(t,x)=\beta_x(t,x)dt + \gamma_x(t,x).dW_t$.
Since $x\mapsto U(t,x)$ is strictly concave and increasing, the marginal utility $ U_x(t,.)
$ is strictly positive and decreasing from $\infty$ to $0$. Consequently $x\mapsto 
U_x(t,x)$ has an inverse  which is the opposite of the marginal conjugate
utility $-\tU_y(t,.)$. By denoting $\mu(t,.):=\beta_x\big(t,-\tU_y(t,.)\big)$
and $\sigma(t,.):=\gamma_x\big(t,-\tU_y(t,.)\big)$, it follows that $U_x$
satisfies the following SDE, with initial condition $U_x(0,x)=u_x(x)$
\begin{equation*}
 dU_x(t,x)=\mu\big(t,U_x(t,x)\big)dt + \sigma\big(t,U_x(t,x)\big).dW_t
\end{equation*}
The change of initial condition from $u_x(x)$ into $z$ yields to the definition
of the process $Z_t(z)=U_x(t,-\tu_y(z))$.
The converse implication is obvious.
\end{proof}
 It remains to give sufficient conditions on the random coefficients
$(\mu,\sigma)$ ensuring the existence of a strong, monotonic solution of SDE
\eqref{eq:ItoSDE}. We briefly recall some classical results on SDEs, useful for
our study; more details and additional results are provided in Section
\ref{subsec:regularityKunita}. An easy to read presentation of SDEs with stochastic coefficients is given in Protter
\cite{Protter}; for more exhaustive study, see Kunita \cite{Kunita:01}.

\begin{theo}\label{theoSDE}
Let us consider the one-dimensional stochastic differential equation, SDE$(\mu, \sigma)$
\begin{equation}\label{ZZZ}
dZ_t=\mu(t,Z_t)dt+\sigma(t,Z_t)dW_t,
\end{equation}
We assume that the $\R$-valued  drift coefficient $\mu=\{\mu(t, x); t \geq 0, x
\geq 0\}$ and the $\R^d$- valued diffusion coefficent $\sigma=\{\sigma(t,x);t
\geq 0, x \geq 0\}$  are
Lipschitz random fields, with random Lipschitz bounds $C_t$ and $K_t$, such that a.s 
$\int_0^TC_tdt<+\infty$ and  $\int_0^TK^2_tdt<+\infty$ for any $T$. In other
words, 
for any $\omega$ outside of a negligible set $N$, for any $t, x, y$,
\begin{equation}
\left\{
\begin{array}{clll}\label{eq: Lipschitzcondition}
|\mu(t,x,\omega)-\mu(t,y,\omega)|\le C_t(\omega)|x-y|,&\quad \mu(t,0)\equiv 0,\\
\|\sigma(t,x,\omega)-\sigma(t,y,\omega)\|\le K_t(\omega)|x-y|,
&\quad  \sigma(t,0)\equiv 0.
 \end{array}
 \right.
 \end{equation}
 
\noindent
\rmi Then, for any $z\in \R_+$ there exists a unique strong solution $\mathbf Z^z $, also called global solution, of
the SDE$(\mu, \sigma)$ \eqref{ZZZ} $(Z^z_0=z)$. Moreover, almost surely, the family of maps $z\mapsto Z_t^z(\omega),~t\ge0$
is continuous
and strictly increasing. \\
\rmii 
 The range of the map $z\mapsto Z(.,z)$ is $]0,+\infty[$ and $Z(.,z)$ is
integrable near to $0$ and to $\infty$. More precisely,
 
 \rma When $z$ goes to $\infty$, almost surely, for any $\eps\in(0,1)$,
uniformly on $[0,T]$ , 
\begin{equation}\label{limit1}
\lim_{z\rightarrow +\infty}\Big(\sup_{0\le t\le
T}\frac{Z(t,z)}{z^{1+\eps}}\Big)=0  \text{ and } \lim_{z\rightarrow
+\infty}\Big(\sup_{0\le t\le T}\frac{Z(t,z)}{z^{\eps}}\Big)=+\infty, \>\text{for
any }T
\end{equation}

 \rmb When $z$ goes to $0$, for any $\eps\in(0,1)$,
 
\begin{equation}\label{limit2}
 \lim_{z\rightarrow 0}\Big(\sup_{0\le t\le T}\frac{Z(t,z)}{z^{\eps}}\Big)=0  
\text{ and } 
\lim_{z\rightarrow 0}\Big(\sup_{0\le t\le
T}\frac{Z(t,z)}{z^{1+\eps}}\Big)=+\infty, \>\>\text{ for all }T
\end{equation}
\end{theo}

\noindent{\bf Comment:} \rmi A constant Lipschitz bound $C$ corresponds to the
classical framework of Lipschitz SDE,  and the assertion (i) is well-known. \\
\rmii The
asymptotic behavior (ii) a) (near to infinity) is less known except in the domain 
 of stochastic flows, where several works and improvements are dealing with this behavior near
of infinity, but unlike Kunita \cite{Kunita:01} and Salah-Eldin \& al \cite{Salah}, Imkeller \& al
\cite{Imkeller} and Zongxia \cite{Zongxia02} consider only  the case of SDEs with deterministic coefficients .
\\
 \rmiii The notion of "global solution" expresses that the solution $ (Z^z_t )$  exists for all $ t \ge0
$. Under weaker assumptions, the solution may be defined  only up to a finite lifetime $ \zeta(z)$. More details will be
given in the next section.\\[3mm]
Sufficient conditions on local characteristics of an  It\^o random
field to be a progressive utility may be exhibited.
\begin{cor} \label{eq:Lipschitzcarac} Assume the framework of Theorem
\ref{UxSDE}. \\
\rmi If there exist random Lipschitz bounds $C_t$ and $K_t$ with
$\int_0^TC_tdt<+\infty$ and  $\int_0^TK^2_tdt<+\infty$ for any $T$, such that
$a.s$, for any $x,x'>0$
\begin{equation}
\left\{
\begin{array}{clll}\label{eq: gamma condition}
|\beta_x(t,x)-\beta_x(t,x')|\leq C_t|U_x(t,x)-U_x(t,x')|, &\lim_{x\rightarrow
\infty}\beta_x(t,x)\equiv 0, \\
 |\gamma_x(t,x)-\gamma_x(t,x')\|\leq K_t|U_x(t,x)-U_x(t,x')|,
&\lim_{x\rightarrow \infty}\gamma_x(t,x)\equiv 0.\\
 \end{array}
 \right.
 \end{equation}
then $\bf U$ is a progressive utility.\\
\rmii Moreover if $(\beta, \gamma)$ are $\cal C^2$-regular random fields,
Condition \eqref{eq: gamma condition}  is equivalent to

\begin{equation}
\left\{
\begin{array}{clll}
|\beta_x(t,x)\leq C_t \,|U_x(t,x)|, &        \|\gamma_x(t,x)\|\leq
K_t\,|U_x(t,x)|             \\
   |\beta_{xx}(t,x)|\leq C^1_t \,|U_{xx}(t,x)|, &\|\gamma_{xx}(t,x)\|\leq K^1_t
\,|U_{xx}(t,x)|\\
 \end{array}
 \right.
 \end{equation}
Then $\mu(t,x)$ and $\sigma(t,x) $ have linear growth with random bounds $C_t$ and $K_t$ respectively, and $\mu_x$ and
$\sigma_x$ are spatially  bounded by $C^1_t$ and $K^1_t$
\end{cor}
\begin{proof} The condition \eqref{eq: gamma condition} is equivalent to the
Lipschitz condition in Theorem \ref{theoSDE} applied to the coefficients $\mu$
and $\sigma$ in Equation \eqref{eq:ItoSDE}.
\end{proof}
\noindent {\bf Remark:} 
In the following, we call  stochastic
flow
any continuous and strictly monotonic solution ${\bf Z}$  of SDEs.  Contrary to
the classical theory, herein we are interested only 
in the process $Z_t\stackrel{def}{=} Z_{s=0,t}$  starting at time $0$, and  its inverse
$(Z_t)^{-1}(y)=Z_{s=0,t}^{-1}(y):=Z_{t,s=0}$, both considered in the forward point of view.
We don't use the general flow associated with the SDE, defined as the family 
$Z_{s,t}(z)$, solution of the equation starting from $z$ at time $s$, and the
flow stability, that is for $r<s<t$, $Z_{s,t}(Z_{r,s}(z))=Z_{r,t}(z)$, except in the paragraph concerning a pathwise dynamic
programming principle (p.\pageref{par:pdpp}).

 From Theorem \ref{UxSDE}, a large part of our study returns to study strictly
monotonic one dimensional random fields and their inverses. The inverse flow of
$\bfU_x$, $(-\bf \tU_y$) plays also a major role in the study of the conjugate
$\bf \tU$ of a  progressive utility $\bf U$.

\subsection{Dynamics of  Convex Conjugate Progressive Utility}\label{DynConjPU}
 The study of the convex conjugate $\bf \tU$ of a progressive utility 
 $\bf U$ is based on the well-known identity (Definition \ref{defCPUF})
$\tU(t,y)=U(t,-\tU_y(t,y))+y \tU_y(t,y)$, and request to know the dynamics  of the ${\mathcal C}^2$-random
field $ U (t, x) $ along the random process $ - \tU_y (t, y) $. For this, an 
extension of the classical It\^o's formula, known as  It\^o-Ventzel's formula
is needed.
We refer to Ventzel \cite{Ventzel} and Kunita \cite{Kunita:01} (Theorem $3.3.1$)  for different variants of this formula and
theirs proofs.

\begin{theo}[It\^o-Ventzel's Formula Weak]\label{IVF}
Consider a ${ \K}^2_{loc}$-It\^o semimartingale $\bfF$ with local characteristics $(\phi,\psi)$.
 For any continuous  It\^o semimartingale $X$, 
$F(.,X_.)$ is a continuous It\^o semimartingale, 
\begin{eqnarray}\label{ItoGen}
&&F(t,X_t)=F(0,X_0)+\int_0^t \phi(s,X_s)ds+\int_0^t\psi(s,X_s).dW_s \\
&+&\int_0^t F_x(s,X_s)dX_s +\frac{1}{2}\int_0^t F_{xx}(s,X_s)\langle dX_s\rangle
+ 
\int_0^t   \langle dF_x(s,x), dX_s \rangle|_{x=X_s} \nonumber
\end{eqnarray}
\end{theo}
\vspace{2mm}

\noindent
{\bf Comment}  The first line of the right hand side of the equation corresponds
to the dynamics of the process $(F(t,x))_{t\ge 0}$ taken on $(X_t)_{t\ge 0}$, 
when in the second line, the first two terms come from the classical It\^o's
formula. The last term  represents the quadratic covariation between $dF_x(t,x)$
and $dX_t$, at $x=X_t$, which can be written  as $\psi_x(t,X_t).\sigma^X_t dt$
when the diffusion coefficient of $X$ is the vector $\sigma^X_t$.

It\^o-Ventzel's formula and monotonic change of variable  will help us to establish
the relationship between  local characteristics of $\bfU$ and $\bf \overline U$.

\begin{theo}\label{ResPrA}
Let $\bfU$  a progressive utility  and $\bf \tU$ its progressive convex conjugate
utility assumed to be ${\K}^2_{loc}$-It\^o semimartingales with local characteristics
$(\beta, \gamma)$ and $(\tilde\beta, \tilde\gamma)$.\\
\rmi The dynamics of  $\bf \tU$ is driven by the non linear second order 
SPDE,
\begin{equation}\label{eq:conjugate}
d\tU(t,y)=\gamma(t,-\tU_y(t,y)).dW_t+\beta(t,-\tU_y(t,y))dt +\frac{1}{2}
\tU_{yy}(t,y)\|\gamma_x\big(t,-\tU_y(t,y)\big)\|^2\,dt
\end{equation}
\rmii The local characteristics of the random field $\bf \tU_y$ are given by, 
\begin{equation}\label{eq:dualcharact}
\left\{
\begin{array}{clllll}
{\tilde \gamma}_y(t,y)&=-\gamma_x(t,-\tU_y(t,y))\tU_{yy}(t,y)\\[1mm]
{\tilde \beta}_y(t,y)&=-\beta_x(t,-\tU_y(t,y))\tU_{yy}(t,y)+
\displaystyle{\frac{1}{2}\partial_y(\frac{\|{\tilde
\gamma}_y(t,y)\|^2}{\,\tU_{yy}(t,y)})}\\[3mm]
\end{array}
\right.
\end{equation}
\rmiii Let $(\mu, \sigma)$ be the random coefficients of the SDE associated with
$\bf U_x$ and 
$\widehat L^{\sigma,\mu}(f)(t,y)=\demi\partial_y(\|\sigma(t,y)\|^2
\partial_yf(t,y))-\mu(t,y) \partial_y f(t,y)$  the adjoint elliptic operator in divergence form associated with
$(\mu,\sigma)$.
Then the marginal conjugate utility $\bf \tU_y$ is a monotonic solution of
the  Stochastic Partial Differential Equation
 with initial condition $\tu_y(y)$,
\begin{eqnarray}\label{AdjointSPDEA}
d\tU_y(t,y)= -\partial_y(\tU_y)(t,y)\sigma(t,y).dW_t+\wL_{t,y}^{\sigma,\mu}(\tU_y) dt
\end{eqnarray}
{\rm Observe that the derivability of the local characteristics 
$(\tilde\beta, \tilde\gamma)$ of $\bf \tU$ requires the existence of a third
derivative for $\bf U$.}
\end{theo}
\begin{proof}  At first, the assumption "$U$ and $\tU$ are ${ \cal
K}^2_{loc}$-semimartingales" implies, from Kunita \cite{Kunita:01} Theorem \ref{DRules} below,  that $(\beta,\gamma)$ and 
$(\tilde\beta, \tilde\gamma)$ belong to the class $ \in {\K}^1_{loc}\times\tilde{\K}^1_{loc}$. Moreover  $U_x$ and
$\tU_y$ are It\^o
semimartingales with local characteristics $(\beta_x,\gamma_x)$ and $(\tilde\beta, \tilde\gamma)$.\\
 Let now  apply the It\^o-Ventzel formula to the
regular random field
$F(t,x)=U(t,x)-y\,x$ and to the semimartingale $X_t=-\tU_y(t,y)$. The following identities related to the change of variable
will be 
useful, 
$F(t,-\tU_y(t,y))=\tU(t,y)$, $U_{xx}(t,-\tU_y(t,y))=-1/\tU_{yy}(t,y)$.\\[1mm] 
\rmi a) Observe that $F_x(t,-\tU_y(t,y))=U_x(-\tU_y(t,y))-y\equiv 0$, so that
the term in $dX$ disappears. Therefore, by It\^o-Ventzel's formula,
the volatility random field ${\tilde \gamma}$ of $\bf \tU$ is ${\tilde
\gamma}(t,y)=\gamma(t,-\tU_y(t,y))$, and its derivative ${\tilde
\gamma}_y(t,y)=-\gamma_x(t,-\tU_y(t,y))\tU_{yy}(t,y)$ is by 	assumption the
volatility characteristic of $\bf \tU_y$. Hence the covariation term is driven
by $\langle dF_x(t,x),-d\tU_y(t,y)\rangle=-\langle \gamma_x(t,x).{\tilde
\gamma}_y(t,y)\rangle dt$.\\
\rmi b)
The  It\^o-Ventzel formula is then reduced to,
\begin{eqnarray*}
d \tU(t,y)&-&\beta(t,-\tU_y(t,y))dt-\gamma(t,-\tU_y(t,y)).dW_t \\
&&= \frac{1}{2}U_{xx}\big(t, -\tU_y(t,y)\big)\langle d\tU_y(t,y)\rangle -\langle
\gamma_x(t,-\tU_y(t,y)).{\tilde \gamma}_y(t,y)\rangle dt\\
&& =\frac{1}{2}U_{xx}(t, -\tU_y(t,y))\|{\tilde \gamma}_y(t,y)\|^2 dt-U_{xx}(t,
-\tU_y(t,y))\|{\tilde \gamma}_y(t,y)\|^2dt\\
&& =-\frac{1}{2}U_{xx}(t, -\tU_y(t,y))\|{\tilde \gamma}_y(t,y)\|^2
dt=\frac{1}{2} \tU_{yy}(t,y)\|\gamma_x\big(t,-\tU_y(t,y)\big)\|^2\
\end{eqnarray*}
\rmii The dynamics of $\bf \tU_y$ is obtained
(by assumption and Theorem \ref{DRules}) by differentiating term by term  in the previous equation.
The use of coefficients $\sigma(t,y)=\gamma_x\big(t,-\tU_y(t,y)\big)$ and
$\mu(t,y)=\beta_x\big(t,-\tU_y(t,y)\big)$ of the SDE associated with $\bf U_x$
allows us to express $\bf \tU_y$ as the solution of a SPDE driven by the
adjoint second order operator in $y$, 
$\wL^{\sigma,\mu}_{t,y}=\demi\partial_y(\|\sigma(t,y)\|^2 \partial_y)-\mu(t,y)
\partial_y$.
\begin{eqnarray*}
d
\tU_y(t,y)&=&-\tU_{yy}(t,y)[\mu(t,y)dt+\sigma(t,y).dW_t]+\partial_y(\frac{1}{2}
\tU_{yy}(t,y)\|\sigma(t,y)\|^2)dt\\
&=& -\partial_y\tU_{y}(t,y)\sigma(t,y).dW_t+\wL^{\sigma,\mu}(\tU_y)(t,y)dt
\end{eqnarray*}
The proof is achieved.\\[-5mm]      
\end{proof} 
\noindent {\bf Remark}
Obviously, we are interested in the properties of the SDE$(\tilde{\mu},\tilde{\sigma})$
 associated with the monotonic random field
 $\bf \tU_y$ when $(\beta, \gamma)$ are $\cal C^2$-regular
random fields.  Given that $\tilde
\sigma(t,-z)=\displaystyle{\frac{\gamma_x(t,z)}{U_{xx}(t,z)}}$ and $\tilde
\mu(t,-z)=
\frac{1}{U_{xx}(t,z)}
\Big(\beta_x(t,z)-\frac{1}{2}\partial_x\big(\frac{\|{\gamma}_x(t,z)\|^2}{\,U_{xx
}(t,z)}\big)\Big)$, it is clear that these coefficients are not globally
Lipschitz and the previous results (Theorem \ref{theoSDE}) cannot be applied directly. So, we report
the study of this SDE in Section \ref{subsec:regularityKunita}, Theorem \ref{SDEInv}, after introducing some
additional tools.

%
\section{Regular Random Fields and  Stochastic Differential Equations}
\label{subsec:regularityKunita}
There are several difficulties in the definition of semimartingales $F(t,x)$
depending on a parameter, as explained in the books of Kunita \cite{Kunita:01} and 
Carmona \& Nualart. \cite{Carmona}, and their local characteristics
$(\phi^F,\psi^F)$ (in short $(\phi,\psi)$  if there is no possible
confusion). \\
\rmi The first one is relative to the existence of continuous version of the
random field  $\bf F$; according to Kunita \cite{Kunita:01} (Theorem 3.1.2
p.75), this property is true when the local characteristics $(\phi,\psi)$ are
locally $\delta$-H\"{o}lder for some $\delta>0$ by Kolmogorov's criterion.\\
\rmii The second one is relative to differential properties:  even if the random
field $\bf F$ and its local characteristics $(\phi, \psi)$ are
differentiable, it is not enough (as is shown in  H. Kunita \cite{Kunita:01}),
to get that the local characteristics of the derivative random field $\bfF_x$ 
are  $(\phi_x, \psi_x)$. \\ 
\rmiii Based on the study of Section \ref{SPU}, we also need under which assumptions on
the coefficients of a SDE, the solution is a regular monotonic random field
with respect to its initial condition.

 We start with a more precise definition of  regular random field spaces than in Section \ref{SPU} by introducing H\"{o}lder
properties. The motivation is find in the Kolmogorov's continuity criterion (\cite{Kunita:01}Theorem 1.4.1).
 \begin{theo} Let $X(x), x\in \mathbb D$ a random field with values in a Banach space $B$, where $\mathbb D$ is a domain in
$\R^d$. Assume that there exist positive  constants $\gamma$, $C$ and $\alpha_i, i=1...d$ with $\sum_{i=1}^d
\alpha_i^{-1}<1$, satisfying,
 \begin{equation}\label{kolmogorov}
 \E[\|X(x)-X(y)\|^\gamma]\leq C\big( \sum_{i=1}^d\|x_i-y_i\|^{\alpha_i}\big), \quad \mbox{ for any}\> x,y\in  \mathbb D
 \end{equation}
 \rmi Then $X(x)$ has a continuous modification $\widetilde X(x)$.\\
 \rmii Let $0<\beta_i\leq A\alpha_i, i=1...d$ arbitrary positive numbers where $A=(\alpha_0-d)/\gamma \alpha_0)$, and
$\alpha_0^{-1}d=\sum_{i=1}^d \alpha_i^{-1}$. Then, for any hypercube $\mathbb H$ there exists a positive random variable 
$K(\omega)$ with $\E[K^\gamma]<\infty$ such that $ \|X(x)-X(y)\|(\omega) \leq K(\omega)\big(
\sum_{i=1}^d\|x_i-y_i\|^{\beta_i}\big) \quad \mbox{ for any}\> x,y\in  \mathbb H, a.s.$.\\[-5 mm]
  \end{theo}
Recently, new criteria based on different norms, or chaining methods have been developed to weaken the criterion. A
interesting survey may be find in Scheutzow \cite{Scheutzow}.

\subsection{ Regular Random Fields Spaces}\label{sub:RSRF}
We introduce a family of Sobolev type random semi-norms to control
locally or globally the growth of the random field  and its
derivatives.  \\[-8mm]
\paragraph{Norms definition}\label{normdef} Let $\phi$ be a continuous $\R^k$-valued progressive random field and let $m$ be
a non-negative integer, and $\delta$ a number in $(0,1]$ . We need to control the asymptotic behavior in $0$ and $\infty$ of
$\phi$, and the regularity of its H\"{o}lder derivatives when there exist. More precisely, 
let $\phi$ be in the class ${\mathcal C}^{m, \delta}(]0,+\infty[)$, i.e.
$(m, \delta)$-times continuously differentiable in  $x$ for any $t$, a.s. \\
\rmi For any  subset  $K\subset ]0,+\infty[$, we define the family of
 random (H\"{o}lder)  $K$-semi-norms 
 \begin{equation}\label{holdernorm}
 \left\{
  \begin{array}{cll}
&\|\phi\|_{m:K}(t,\omega)=\sup_{\substack{x\in
K}}\frac{\|\phi(t,x,\omega)\|}{x}+\sum_{1\le j\le m}\sup_{\substack{x\in
K}}\|\partial_x^{j}\phi(t,x,\omega)\|\\[2mm]
& \|\phi\|_{m,\delta:K}(t,\omega)=\|\phi\|_{m:K}(t,\omega)+\displaystyle{\sup_{\substack{x,
y\in
K}}\frac{\|\partial_x^{m}\phi(t,x,\omega)-\partial_x^{m}\phi(t,y,\omega)\|}{
|x-y|^\delta}}.
\end{array}
\right.
 \end{equation}
The case $(m=0, \delta=1)$
corresponds to the local version of the Lipschitz case used in Section 1.
When $K$ is  all the domain $ ]0,+\infty[$, we simply write
$\|.\|_{m}(t,\omega)$, or $\|.\|_{m,\delta}(t,\omega)$.\\
\rmii The first term of these random  semi-norms differs slightly from
the definition of
Kunita semi-norms (Equations (1) and (2) p.$ 72$) because  instead of dividing
by $1+|x|$ we divide by $x$ on the first terms. This does not change   Kunita's
results, but  allows us to obtain reasonable behavior  in the neighborhood of
$x= 0$  (Equation \eqref{limit2}, Theorem \ref{theoSDE}) in addition to the
traditional results in the neighborhood of $x= \infty$. \\[-8mm]
\paragraph{Different spaces of regular random fields}\label{par:spacesregular}
\noindent 
The previous semi-norms are related to the spatial parameter. As in Definition \ref{Kregularity}, we add the temporal
dimension in assuming
 these semi-norms (or the square of the semi-norm) to be
integrable in time with respect to the Lebesgue measure on $ [0, T] $ for all $
T $. Then, as in Lebesgue's Theorem, we can differentiate, pass to the limit, commute
limit and integral for the random fields.  Calligraphic notation recalls that these semi-norms are random. \\
\rmi
 ${\K}^{m,\delta}_{loc}$ (resp. $\overline{{\K}}^{m,\delta}_{loc}$) denotes the set of all  $\Cc^{m,\delta}$-random
fields
such that  for any compact $K\subset]0,+\infty[$, and any $T$,
$\int_0^T\|\phi\|_{m,\delta:K}(t,\omega)<\infty$, (resp.$ \int_0^T\|\psi\|^{
2}_{m,\delta:K}(t,\omega)dt<\infty$ ). \\
\rmii  When these different norms are well-defined on the whole space $]0,+\infty[$, the
derivatives (up to a certain order) are bounded in the spatial parameter, with
integrable  (resp.square integrable) in time random bound. In this case, we use the notations
${\K}^{m}_b,\overline{{\K}}^{m}_{b}$ or ${\K}^{m,\delta}_b,\overline{{\K}}^{m,\delta}_{b}$.

\subsection{Differentiability of It\^o random fields}\label{FRSDE}
We shall discuss the regularity of a It\^o semimartingale random field
$$\bf{F}(t,x)=F(0,x)+\int_0^t\phi(s,x)ds+\int_0^t\psi(s,x).dW_s$$
in connection with the regularity of its local characteristics $(\phi,\psi)$. As in Section 1, by convention, an It\^o random
field $\bfF$ is said to be a ${ \K}^{m,\delta}_{loc}$-{\bf semimartingale}, whenever $F(0,x)$ is of class $\Cc^{m,\delta}$,
$B^F(t,x)=\int_0^t\phi(s,x)ds$ is of class ${ \K}^{m,\delta}_{loc}$, and $ M^F(t,x)=\int_0^t\psi(s,x).dW_s$ is of class ${
\oK}^{m,\delta}_{loc}$. The reference to ${ \K}^{m,\delta}$ recall that $\bfF$ is a random field.
\footnote{
Note, in the case
of general semimartingales $F$ studied by Kunita, the notion of local
characteristics makes reference to the triplet $(\phi(t,x),a(t,x,y),A_t)$ for
$x,y >0$, where the joint quadratic variation of $F(t,x),~F(t,y)$ satisfy the
relation $<F(t,x),F(t,y)>=\int_0^ta(s,x,y)dA_s$ with the process $A$ is 
$\F-$adapted increasing and continuous. Here, as there is a finite number of 
Brownian motions, $A_t=t$ and $a(t,x,y):=\psi(t,x).\psi(t,y)$ where, as
before,  the  "." denote the inner scalar product. From this, all  assumptions
on the matrix $"a"$ (in Kunita) are adapted to the present framework and replaced by our
equivalent hypothesis on the  diffusion characteristic vector $ \psi $.}
As in Kunita \cite{Kunita:01}, we are concerned  both by the regularity of $F$ 
from the regularity of  its local characteristics $(\phi,\psi)$ (Theorem 3.1.2 ) and by the regularity of  $(\phi,\psi)$ from
that of $\bf{F}(t,x)$ (Theorem
3.1.3). To be concise, we also give a sufficient conditions  (\cite{Kunita:01} Theorem $3.3.3$) under
which we can differentiate term by term the dynamics of an It\^o random field.  This property is used in order to apply
It\^o-Ventzel's formula. 
\begin{theo}[Differential Rules]\label{DRules}
Let  $\bfF$  be an It\^o semimartingale random field with local characteristics  $(\phi,\psi)$, 
 $F(t,x)=F(0,x)+\int_0^t\phi(s,x)ds+\int_0^t\psi(s,x).dW_s$\\[1mm]
\rmi If $\bfF$ is a ${\K}^{m,\delta}_{loc}$-semimartingale for some $m\ge0,~\delta \in (0,1]$,  its local characteristics 
$(\phi,\psi)$ are of class ${\K}^{m,\eps}_{loc}\times \overline{{\K}}^{m,\eps}_{loc}$ for any $\eps<\delta$.\\[1mm]
 \rmii Conversely, if the local characteristics $\mathbf{ (\phi,\psi)}$ are of class  ${\K}^{m,\delta}_{loc}\times
\overline{{\K}}^{m,\delta}_{loc}$,
then $\bfF$ is a ${\K}^{m,\eps}_{loc}$-semimartingale for any $\eps<\delta$.\\
\rmiii In any cases, for $m\ge 1, \delta \in (0,1]$, the derivative random field $\bfF_x$ is an It\^o random field with local
characteristics $(\phi_x,\psi_x)$.\\
\rmiv Moreover, if  $\bfF$  is a ${\K}^{1,\delta}_{loc}\cap {\Cc}^{2}$-semimartingale,  for any It\^o process
$X$, $F(.,X_.)$ is a continuous It\^o semimartingale satisfying the It\^o-Ventzel formula  \eqref{ItoGen}.
\end{theo}
 \noindent As previously mentioned, we also need results on the existence and the regularity of
one dimensional random fields which are also solutions of stochastic differential equations (SDE).
The spatial parameter in this case corresponds to the initial condition.  Such random fields are also called stochastic flows
and are the main subject (in the multidimensional case) of the Kunita's book \cite{Kunita:01}.\\
The question is now to make assumptions on the coefficients in place of local
characteristics. An example was given in Section \ref{SPU}, Theorem \ref{theoSDE}, where the
global Lipschitz regularity of  coefficients is used in proving the existence of monotonic solution.
Next proposition, (\cite{Kunita:01} Theorem $4.6.5$), frequently used in the sequel, extends these results to
differentiability
properties. 
\begin{pp}\label{SDETHEO}
Let  $\mu$ be a real valued process and $\sigma$ be a $d$-dimensional process of the
  the class $ {\K}^{m,\delta}_{\bf b}$ and $\overline{{\K}}^{m,\delta}_{\bf b}$ for some  $ m\ge 1,\delta\in (0,1]$.
Consider the following SDE$(\mu,\sigma),$
\begin{equation}\label{EDSZ}
dX_t=\mu(t,X_t)dt+\sigma(t,X_t).dW_t,\quad X_0=x
\end{equation}
Then, in addition  to the results of Theorem \ref{theoSDE} $\big((\mu, \sigma)\in  {\K}^{0,1}_{\bf b}\times
\overline{{\K}}^{0,1}_{\bf b}\big)$ on
the existence and uniqueness of global monotonic solution, we
have:\\
\rmi The unique
solution ${\bf X}=(X_t^x, x>0)$  is a $ \K^{m,\eps}_{loc}$ semimartingale
for any $\eps<\delta$. The
inverse ${\bf X^{-1}}$ of $\bf X$ is also of class ${\cal C}^{m}$. \\[1mm]
 \rmii  ${\bf X}$ is strictly increasing in $x$, and its derivative $\bf X_x$
is solution of a linear equation, with spatially bounded  stochastic
parameters
$(\mu_x(t,X_t^x),\sigma_x(t,X_t^x))$ given by
\begin{equation}\label{Z_zF}
 dX_x(t,x)=X_x(t,x)\big[\mu_x(t,X_t^x)dt +\sigma_x(t,X_t^x).dW_t\big], \> X_x(0,x)=1
\end{equation} 
$\bf{1/X_x}$ is solution of the same kind of linear equation 
 \begin{equation}\label{Z_zF'}
 dZ_t=Z_t\big[(-\mu_x(t,X_t^x)+||\sigma_x(t,X_t^x)||^2)dt -\sigma_x(t,X_t^x).dW_t\big]
\end{equation} 
Then, $\bf X_x$ and $\bf 1/X_x$ are $ {\K}^{m-1,\eps}_{loc}$-semimartingales.\\
\noindent
\rmiii The local characteristics of ${\bf X}$, $\lambda^X(t,x)=\mu(t,X^x_t)$ and
$\theta^X(t,x)=\sigma(t,X^x_t)$ have only local properties and belong to ${\K}^{m,\eps}_{loc}\times
\overline{{\K}}^{m,\eps}_{loc}$ for any $\eps<\delta$. 
\end{pp}
\begin{req}
 {\rm This technical result  shows clearly the interest  of using  H\"{o}lder property: the solution is fractionally less
regular than the coefficients (going from $ \delta $ in $ \eps <\delta $).
Otherwise, if we are only interested with  processes of class $ \K^m$ ($m$ integer) without worrying about the H\"{o}lder's
dimension, then we will lose a hole order in the regularity: instead of  solution of class $\K^{m,\eps}_{loc}$, we
will only obtain solution of class  $\K^{m-1}_{loc}$.}
\end{req}
\noindent In the last proposition, we have made global regularity assumptions on the coefficients $(\mu,\sigma)$ in order to
prove that
the SDE \eqref{EDSZ} has a regular, non-exploding (global), strictly monotonic solution $X$. Observe that the solution itself
as random field is only in ${\K}^{m,\eps}_{loc}$ in general. In Section \ref{SPU}, Paragraph \ref{DynConjPU}, we have shown
that even if
the progressive marginal utility satisfies a SDE$(\mu,\sigma)$ with global regularity, the SDE$(\tilde\mu,\tilde\sigma)$
satisfied by the inverse $-\tU_y$ of $U_x$, assumed to be  an It\^o process,
has only local regularity. So, we are also concerned with SDE whose the coefficients are of class
${\K}^{0,1}_{loc}\times\overline{{\K}}^{0,1}_{loc}$ or 
${\K}^{m,\delta}_{loc}\times\overline{{\K}}^{m,\delta}_{loc}, ~m\ge1,~\delta\in ]0,1]$. Then, the SDE can not have global
solution and explosion can occur at  finite time $\zeta(x).$ Nevertheless, several properties are maintained as shown in
\cite{Kunita:01} Theorems $4.7.1$ and $4.7.2$.
\begin{theo}[Local assumptions]\label{SolMax}
Let us consider a SDE  with only  locally Lipschitz coefficients
$\big( (\mu,\sigma)\in {\K}^{0,1}_{\bf loc}\times \overline{{\K}}^{0,1}_{\bf loc}\big)$.\\
 \rmi  For any initial
condition $x$, the SDE  has a unique  maximal solution $(X^x_t)$ up to an
explosion time $\zeta(x)$, that is on $[0,\zeta(x))$,  
\begin{equation}\label{eq: explosive}
 dX_{t}^x=
\ind_{\{t<\zeta(x)\}} \mu(t,X^x_{t})dt
+\ind_{\{t<\zeta(x)\}}\sigma(t,X^x_t)dW_t, \quad a.s.
\end{equation}
\rmii For given $t$, as function of $x$,  $X^x_t$ is defined on its domain
 ${\cal D}_t=\{x: \zeta(x)>t\}$, with finite values on its range ${\cal R}_t(\omega)=\{X_t^x(\omega); x\in {\cal
D}_t(\omega)\}$.
Furthermore, $x\mapsto X^x_t :{\cal D}_t\rightarrow{\cal R}_t $ is a continuous 
strictly monotonic random field with continuous inverse $X_t^{-1}: {\cal R}_t\rightarrow{\cal D}_t $.\\
\rmiii $(X^x_t)$  is a global solution (Theorem \ref{theoSDE}) if and only if the explosion time $\zeta(x)$ is equal to
$\infty$
for all $x\in\R^+=[0,\infty)$ a.s., (or equivalently if for any time $t$, the domain ${\cal D_t}$ 
is the whole space $\R^+$ a.s.). Hence the range ${\cal R}_t$ is also the whole space $\R^+$ for any $t.$ \\
\rmiv Moreover, if $(\mu,\sigma)\in{\K}^{m,\delta}_{loc}\times\overline{{\K}}^{m,\delta}_{loc}$ $(m\ge1,~\delta\in
(0,1])$  $X_t(.)$ is of class $\Cc^{m,\eps}, \eps<\delta$ on ${\cal D}_t$. \\
\rmv  If $(X^x_t)$  is a global solution,  then all  assertions of
Proposition \ref{SDETHEO} hold true.
\end{theo}
\noindent We recall the proof of Kunita to present the truncation method which is the base of many proofs when only local
properties hold true.
\begin{proof}
\rmi Following Kunita, we shall apply the method of truncation on $\R^+$. For each positive integer $N$, take a
$\Cc^\infty$-function
$\psi^N(x)$, $x\in\R^+$
such that $\psi^N(x)=1$ if $x\le N$, $0\le\psi^N(x)\le1$ if $N\le x\le N+1$ and $\psi^N(x)=0$ if  $x> N$. Define
$\mu^N(t,x)= \mu(t,x)\psi^N(x)$ and $\sigma^N(t,x)= \sigma(t,x)\psi^N(x)$ which belong to  ${\cal
L}^{m,\delta}_{b}$ and $\overline{{\K}}^{m,\delta}_{b}.$ Therefore, the SDE$(\mu^N,\sigma^N)$
has a unique solution $X^N$ of class $\K^{m, \epsilon}$ $(\epsilon<\delta)$ satisfying assertions of Proposition
\ref{SDETHEO}. 
Let us introduce the family of stopping times for each $x$, $\zeta_N(x):=\inf\{t: X^N(t,x)\ge N \}$.
 By uniqueness, for $M<N$ the solutions $X^M_.(x)$ and $ X^N_.(x)$ coincide on $[0,\zeta_M(x))$. Then, the family
$\zeta_N(x)$ is increasing with limit $\zeta(x):=\lim_{N \rightarrow\infty}\zeta_N(x)$.
We can define $X_t(x),~t<\zeta(x)$ by  
$X_t(x)=X^N_t(x)$ if $t<\zeta_N(x)$. It is a maximal solution of SDE$(\mu,\sigma)$ starting at $x$ at
$t=0$.\\
 \rmii We shall prove that the above  $X_.(x)$ is a continuous strictly increasing flow. Take any sample $\omega$ such that
$X^N_t(\omega,.)$ defines a continuous map strictly increasing in $x$ for any $N$. Probability of the set of all such
samples
is 1. We first note that each $\zeta_N(\omega,x)$ is lower semicontinuous in $x$ i.e. ${\cal
D}^N_t(\omega)=\{x:\zeta_N(\omega,x)>t\}$ is open for
any $t>0$, since on this set $X^N_s(\omega,x)<N$ for all $s\le t$. 
Hence, the same inequality holds for
any $x'$ in the neighborhood of $x$. Now since $\zeta(\omega,x)$ is the upper limit of $\zeta_N(\omega,x)$ it is also lower
semicontinuous. Set ${\cal D}_t(\omega)=\cup_N {\cal D}^N_t(\omega)$.
Since $X(t,x)=X^N(t,x)$ holds on ${\cal D}^N_t(\omega)$, the map $X(t,x):{\cal
D}_t(\omega)\rightarrow \R^+$ is continuous strictly increasing with continuous inverse.\\
 \rmiii The  differentiability properties are showed similarly by truncation techniques, using derivability of
processes $X^N$.
\end{proof}
\noindent {\bf Comment} Recently, several papers address the question: under which minimal assumptions on the  SDE's 
coefficients, the
solution is non-explosive? For example, when the coefficients $\mu$ and $\sigma$ are a deterministic functions independent of
the time, the property
holds true under global log-Lipschitz type conditions  as it is showed  in  Zongxia \cite{Zongxia01} and Fang
\cite{Fang}. But these
new results and many others can not be applied directly to our study because the coefficients of SDE's we are concerned 
 are structurally stochastic.

Local regularity  on SDEs coefficients appears as a kind of minimal assumption  to ensure   the
regularity
of a global solution {\em if there exists}.  Because of its importance in the sequel,  we give a name to this
class of SDE's.
\begin{df}\label{def:EDSSolution}
A SDE$(\mu, \sigma)$ is said to be of class $ \S^{m,\delta}$ if

\rma  the coefficients $(\mu, \sigma)$  are in the spaces  $({\K}^{m,\delta}_{loc},\overline{{\K}}^{m,\delta}_{loc}))$

\rmb  the  maximal solution $X$ is non explosive.\\[1mm]
By Theorem \ref{SolMax}, the unique solution $X$ is strictly monotonic with range $[0,\infty)$ and of class
$\K^{m,\eps}_{loc},~\eps<\delta$.

\end{df}

\noindent Classical examples of $\S^{m,\delta}$ SDEs
 are given by SDE$(\mu, \sigma)$ when  $(\mu, \sigma)$ are in the spaces $({\K}^{m}_{b},\overline{{\K}}^{m}_{b})$, or
even in
$({\K}^{0}_{b},\overline{{\K}}^{0}_{b}) \cap({\K}^{m,\delta}_{loc},\overline{{\K}}^{m,\delta}_{loc})$.
Moreover, in these
last two cases, the asymptotic behavior of the solution is given by Theorem \ref{theoSDE}, Equation \eqref{limit1} and
Equation \eqref{limit2}. 
%
\subsection{Solvable SPDEs via SDEs}\label{DIFRSDE}
\paragraph{Dynamics of  inverse flow of regular SDE solution }
\noindent
 In the utility framework, under the strong assumptions of  Theorem \ref{ResPrA},
we have shown that $\tU_y$, the inverse flow of $(-U_x(t,x))$, is 
solution of a SPDE and a SDE$(\tilde \mu, \tilde \sigma)$ simultaneously. 
We want to relax the a priori assumption that the inverse flow  $\tU_y$ is an It\^o random field. 
So, we proceed differently,  by starting from the maximal solution of the local SDE $(\tilde \mu, \tilde
\sigma)$ with explosion time $\zeta(x)$,  and by verifying that up to $\zeta(x)$ this process is the inverse flow of the 
solution of the SDE$( \mu, 
\sigma)$. By an easy argument base on the uniqueness of non-explosive solution of the SDE$(\mu,\sigma)$, we deduce
$\zeta(x)=\infty,
\>a.s.$. Some regularity on the SDE$(\mu,\sigma)$ is required to conduct calculation and conclude.
%
\begin{theo}[Inverse flow SDE]\label{SDEInv}
Let $(X^x_t)$ be the monotonic solution of a SDE$(\mu,\sigma)$ of  class $\S^{m, \delta},m\ge 3,~\delta\in ]0,1]$,
so that as
random
field
 $(X^x_t)$ and its local characteristics $(\lambda(t,x)=\mu(t,X^x_t))$ and $(\theta(t,x)=\sigma(t,X^x_t)) $ are of class
$\K^{m,\eps}_{loc}$  and  $ {\cal
L}^{m,\eps}_{loc}\times\overline{{\K}}^{m,\eps}_{loc}$ for any $0<\eps<1$.\\
 We are concerned with
the SDE$( \tilde\mu,\tilde\sigma)$ 
\begin{equation}\label{EqFInv}
d\xi_t=-\frac{1}{X_x(t,\xi_t)}\Big[\big(\lambda(t,\xi_t)-\frac{1}{2}
\partial_x\big(\frac{ \|\theta\|^2}{X_x}\big)(t,\xi_t)\big)dt+
\theta(t,\xi_t).dW_t\Big],~ \xi_0=z,
\end{equation}
where 
 $\displaystyle \tilde \sigma(t,z)=-\frac{\theta(t,z)}{X_x(t,z)}$ and 
$\displaystyle \tilde \mu(t,z)=\frac{1}{X_x(t,z)}\Big(\frac{1}{2}
\partial_x\big(\frac{ \|\theta\|^2}{X_x}\big)(t,z)-\lambda(t,z)\Big). $\\[3mm]
 \rmi The SDE$(\tilde \mu,\tilde \sigma)$ is of  class $ \S^{m-2,\eps}$$(0<\eps <\delta)$ and its unique monotonic solution
$\xi^z$ is the  inverse flow $X^{-1}$ of $X$.\\
\rmii Consequently, the inverse $X^{-1}$ of $X$ is a semimartingale and 
belongs to the class ${\K}^{m-2,\eps}_{loc}\cap {\cal C}^{m}.$
\end{theo}

\begin{proof} The proof is in several steps, by first proving the local regularity of the coefficients $(\tilde \mu,\tilde
\sigma )$, and then the existence of a monotonic solution $\xi^z$ up to a explosion time $\zeta(z).$ The main step is then to
prove that $\zeta(z)=\infty \>a.s.$ by showing that locally $\xi^z$ is the inverse $X^{-1}$.\\[1mm]
\rmi 
Since $X\in{\K}^{m,\eps}_{loc}$ and $(\lambda,\theta)\in {\K}^{m,\eps}_{loc}\times\overline{{\K}}^{m,\eps}_{loc}$, by
Proposition \ref{SDETHEO}
$1/X_x\in{\K}^{m-1,\eps}_{loc}$ and so
$(\lambda/X_x,\theta/X_x)\in  {\K}^{m-1,\eps}_{loc}\times\overline{{\K}}^{m-1,\eps}_{loc}$ and 
$ \partial_x(\frac{ \|\theta\|^2}{X_x})\in {\K}^{m-2,\eps}_{loc}. $ Consequently  the coefficients $(\tilde \mu,\tilde
\sigma)$ 
are of class $ {\K}^{m-2,\eps}_{loc}\times \overline{{\K}}^{m-1,\eps}_{loc}.$ \\[1mm]
\rmii  Since the coefficients $(\tilde{\mu},\tilde{\sigma})$ satisfy the assumptions of Theorem \ref{SolMax}, 
the SDE$(\tilde{\mu},\tilde{\sigma})$ has a unique maximal solution $\xi^z$, up to an explosion
time $\zeta(z)$.\\
 \rma We claim that the solution $\xi_t^z$ is the inverse flow $X^{-1}(t,z)$ for
$t\in [0,\zeta(z))$ of the SDE$(\mu, \sigma)$ monotonic solution $X_t^x:=X(t,x)$.
Since by assumption $X$ is of  class $\K^{m, \epsilon}_{loc}$ and its local characteristics $(\lambda,\theta)$ are  of class
$(\lambda,\theta)\in {\K}^{m,\eps}_{loc}\times
\overline{{\K}}^{m,\eps}_{loc}$ $(m\geq 3, \epsilon \in ]0,\delta[),$
we can 
apply It\^o-Ventzel's formula to $X(t,\xi_t^z)$ up to the time $\zeta(z)$.
Then, on $[0,\zeta(z))$, by using the short notation $\xi$ in place of $\xi^z$, \begin{eqnarray*}
 dX(t,\xi_t)&=&\Big[\theta(t,\xi_t)+X_x(t,\xi_t)\frac{
-\theta(t,\xi_t)}{X_x(t,\xi_t)}\Big].dW_t\\
&+&
 \Big[\lambda(t,\xi_t)+X_x(t,
\xi_t)\Big(-\frac{1}{
X_x(t,\xi_t)}\Big(\lambda(t,\xi_t)-\frac{1}{2}\partial_x\big(
\frac{ \|\theta\|^2}{X_x}\big)(t,\xi_t)\Big)\Big)\Big]dt\\
 &+&\Big[
\demi X_{xx}(t,\xi_t)\big(\frac{\|\theta\|^2}{X_x^2}\big)(t,\xi_t)+\big(\theta_x(t,\xi_t)
.\big(\frac { -\theta } {X_x}\big)(t,\xi_t)\Big]dt \\
&=&\big[0\big].dW_t+\Big[\frac{1}{2}\partial_x\big(
\frac{ \|\theta\|^2}{X_x}\big)(t,\xi_t)-\demi \big(\frac{1}{X_x}\partial_x\|\theta\|^2\big)(t, \xi_t))-\demi
\big(\partial_x(\frac{1}{X_x})\|\theta\|^2
\big)(t,\xi_t)\Big]dt\\
&=&0
\end{eqnarray*}
\rmb Then the continuous (in time) process $X(t,\xi^z_t)$ is constant a.s. on $[0,\zeta(z))$.
At time $t= \zeta(z)<\infty$,  $\xi_t^z=\infty$ and $X(t,\infty)=\infty$. On the other hand, by continuity, $X(t,\xi^z_t)=z$
if  $t=
\zeta(z)<\infty$.
To avoid contradiction, necessarily $\zeta(z)=\infty ,\>a.s.$.
\end{proof}
\noindent As in Theorem \ref{ResPrA}, we can also characterize the inverse process in terms of monotonic solution of
non-linear
stochastic partial differential
equation (SPDE). This point of view is well-suited  to the study of consistent dynamic
utilities developed in the next sections. \\[-6mm]
\paragraph{Solvable SPDEs via SDEs}
Now we show how to solve some  non-linear SPDEs via SDEs. A first link  is
obtained by the transformation of  the SDE \eqref{EqFInv} into a SPDE. Recall that conditions of type $\cal S^{m,\delta}$ are
related to the coefficients of the SDE, when conditions of type $\K_{loc}^{m,\delta}$ are related to the
local characteristics of the solution.
\begin{theo}[SPDE point of view]\label{SPDEInv} 
Let us consider a SDE $(\mu,\sigma)$ of class ${\cal S}^{m, \delta}$
 with $m\geq 2, \delta \in (0,1]$, and  its adjoint operator
$\wL^{\sigma,\mu}_{t,z}=\demi\partial_z(\|\sigma(t,z)\|^2 \partial_z)-\mu(t,z) \partial_z$. Denote by  $X$  its unique
solution. \\
\rmi For $m\geq 3$, the inverse flow 
$X^{-1}=\xi^X$ of $X$ is a  strictly monotonic solution of  class ${\K}^{m-2,\delta}_{loc}\cap\cal C^m$ of 
 SPDE$(\widehat L^{\sigma,\mu}, -\sigma \partial_z)$, with initial condition $\xi_0(z)=z$,
\begin{eqnarray}\label{AdjointSPDE}
d\xi(t,z)= -\xi_z(t,z)\sigma(t,z).dW_t+\wL_{t,z}^{\sigma,\mu}(\xi) dt
\end{eqnarray}
\rmii Conversely,  $(m\geq 2)$, let $\xi$ be a ${\K}^{1,\delta}_{loc}\cap{\cal C}^{2}$-regular  solution of  
SPDE$(\wL^{\sigma,\mu},-\sigma \partial_z)$ \eqref{AdjointSPDE}. Then,
$\xi(t,X(t,x))\equiv x$ and $\xi$ is  the strictly monotonic inverse flow  $X^{-1}:=\xi^X$ of $X$.
 Moreover, uniqueness  holds true for the SPDE$(\wL^{\sigma,\mu},-\sigma \partial_z)$ in the class of 
${\K}^{1,\delta}_{loc}\cap \cal C^2$-regular solutions. 
\end{theo}
\noindent {\bf Comment} When the coefficients $(\mu, \sigma)$ of the SDE are non random, in a multi-dimensional case, the
SPDE$(\wL^{\sigma,\mu},-\sigma \partial_z)$ is known as the stochastic transport equation in Zhang \cite{Zhang} or in
Constantin and Iyer \cite{Constantin}.
\begin{proof}
\rmi According to notations of Theorem \ref{SDEInv}, since $\xi^X$ is the inverse of $X$,  we have:
\begin{eqnarray*}
\lambda(t,\xi^X(t,z))&=&\mu\Big(t,X\big(t,\xi^X(t,z)\big)\Big)=\mu(t,z), \qquad \xi^X_z(t,z)=\frac{1}{X_x(t,\xi^X(t,z))}, \\
 \theta(t,\xi^X(t,x))&=&-\frac{\sigma\Big(t,X\big(t,\xi^X(t,z)\big)\Big)}{X_x(t,\xi^X(t,z))}=-\xi^X_z(t,z)\sigma(t,z)
\end{eqnarray*}
which easily leads to $ \tilde{\mu}(t,\xi^X(t,z))= \widehat L_{t,z}^{\sigma,\mu}(\xi^X)(t,z).$ So $\xi^X$ is
solution of the SPDE \eqref{AdjointSPDE}.\\
\rmii Let now turn to the converse implication by starting from a given
monotonic solution  $\xi$ of class $\K^{1,\delta}_{loc}\cap \cal C^2$ of the SPDE: $d\xi(t,z)=
-\xi_z(t,z)\sigma(t,z).dW_t+\wL_{t,z}^{\sigma,\mu}(\xi)
dt$. \\
\rma From Theorem \ref{DRules}, $\xi$ is regular enough to use It\^o-Ventzel's formula with the
 solution  $X(t,x)=X^x_t$ of the SDE$(\mu,\sigma)$ to compute the dynamics of $H(t,x)=\xi(t,X(t,x))$. In the next equation,
we do not recall the parameter $x$.
\begin{eqnarray*} 
 dH_t&=&\big(-\xi_z(t,X_t)\sigma(t,X_t)-\xi_z(t,X_t)\sigma(t,X_t) \big).dW_t\\ \nonumber
 &+&\Big(\wL^{\sigma,\mu}(\xi)+  \demi \xi_{zz}\|\sigma\|^2 +\mu \>\xi_z+\partial_z(-\xi_z\sigma).\sigma\Big)(t,X_t)dt\\
 &=&\Big(\xi_{zz}\|\sigma\|^2+\demi \xi_z(\partial_z\|\sigma\|^2)-\partial_z(\xi_z)\|\sigma\|^2-\demi
\xi_z(\partial_z\|\sigma\|^2)\Big)(t,X_t)dt\\
 &=&0
 \end{eqnarray*}
 The random field $H(t,x)=\xi(t,X(t,x))$ is constant in time and equal to its initial condition $x$. This finishes the proof
that $X$ is the inverse flow of $\xi$.\\
\rmb Since SDE$(\mu,\sigma)$ is of $\cal S^{m,\delta}$, there is only one solution $X$. Then any
"regular" solution 
$\xi$ of the SPDE is the inverse of $X$ and then is unique. 
\end{proof}
\noindent Next result, useful for applications, is a slight extension of the previous one.
It establishes a connection between a more general second order  SPDE and two SDEs. As discussed in
the following, this connection is in the core of the study of consistent dynamic utilities. It is based on the observation
that
if $\xi$ is the inverse of SDE$(\mu^X,\sigma^X)$ monotonic solution $X$ and if $\phi\in \Cc^2$ a regular monotonic
function, the process $X(.,\phi(x))$ satisfies the same SDE$(\mu^X,\sigma^X)$, and so  its inverse 
$\phi^{-1}(\xi_.(z))$ satisfies the same SPDE than $\xi$. The extension, given in the following result describes the dynamics
of compound processes $Y(t,\xi(t,z))$ and  $\xi(t,\overline{X}(t,\overline{x})$ for regular It\^o semimartingales $Y$ and
$\overline{X}$ and show how to solve
the  associated SPDEs.

\begin{theo}\label{pp: composition formula}  Let $X$ be a solution of SDE$(\mu^X,\sigma^X)$ and $\xi$ a 
${\K}^{1,\delta}_{loc}\cap \cal C^2$-regular solution $(\delta>0)$ of the
SPDE$(\widehat L^{X},-\sigma^X \partial_z)$, where $\widehat L_{t,z}^{X}=\widehat
L_{t,z}^{\sigma^X,\mu^X}=\demi\partial_z(\|\sigma^X(t,z)\|^2
\partial_z)-\mu^X(t,z) \partial_z$.\\[1mm] 
\rmi Let $Y$ be a solution of class ${\K}^{1,\delta}_{loc}\cap \cal C^2$ of SDE$(\mu^Y,\sigma^Y)$ with initial condition
$\phi\in \Cc^2$. Then the random field $Y(t,\xi(t,z))=G(t,z)$ evolves as,
\begin{eqnarray}\label{eq: composition formulaA}
 dG(t,z)&=&\sigma^Y(t,G(t,z)).dW_t+\mu^Y(t,G(t,z))dt\nonumber \\
 &-& \partial_zG(t,z)\sigma^X(t,
z)\big[dW_t+\sigma^Y_y(t,G(t,z))dt]+ {\widehat L}_{t,z}^{X}(G)(t,z))dt
\end{eqnarray}
 \rmii Let $\overline{X}$ be a solution of SDE$(\overline{\mu},\overline{\sigma})$
 with initial condition $\psi(\overline{x})$. Denote by
$\Delta{\mu}(t,z):=\mu^X(t,z)-\overline{\mu}(t,\overline{X}_t)$ and $\Delta{\sigma}(t,z):=
\sigma^X(t,z)-\overline{\sigma}(t,\overline{X}_t)$. Then 
the random field  $\xi(t,\overline{X}_t)$ evolves as
 \begin{eqnarray}\label{eq: compcharacteristic}
 d\xi(t,\overline{X}_t)=-\xi_z(t,\overline{X}_t)\big[\big(\sigma^X(t,\overline{X}_t)-\overline{\sigma}(t,\overline{X}
_t)\big).(dW_t-\overline{\sigma}_{\overline{x}}(t,\overline{X}_t)dt)\big]
 +{\widehat L}^{\Delta}(\xi)(t,\overline{X}_t).
 \end{eqnarray}
  \rmiii {\bf Solvable SPDE:} Let $G$ be a  solution of class ${\K}^{1,\delta}_{loc}\cap \cal C^2$-regular of the
SPDE \eqref{eq: composition formulaA}; then
the process $G(t,\overline{X}_t(\overline{x}))$ with initial condition $G(0,\psi(\overline{x}))$ evolves as 
 \begin{eqnarray}\label{eq: compcharacteristicA}
 dG(t,\overline{X}_t)&=&\sigma^Y(t,G(t,\overline{X}_t)).dW_t+\mu^Y(t,G(t,\overline{X}_t))dt\\
 &-&G_z(t,\overline{X}_t)\big(\sigma^X-\overline{\sigma})(t,\overline{X}_t).\big[dW_t+(\sigma^Y_y(G)-\overline{\sigma}_{
\overline{x}})(t,\overline{X}_t)dt\big]
+ {\widehat L}^{\Delta}(G)(t,\overline{X}_t)dt\nonumber
\end{eqnarray}
\rmiv
 In particular, $G(t,X_t(\phi(y))$ is a solution of the SDE$(\mu^Y,\sigma^Y)$ with initial condition $\phi(y)$. If
uniqueness holds true for this equation, then
$G(t,z)=Y_t(t,\xi(t,z))$ and uniqueness also holds true for the SPDE \eqref{eq: composition formulaA}.
  \end{theo}
\noindent  Note the different nature of assumptions (which  may be equivalent) in the
assertions of this theorem. In $\rm (i)$, we assume that the
coefficients are regular enough such that $Y$ satisfies the It\^o-Ventzel
assumptions and such that  the inverse $ \xi $ of $ X $ is an It\^o semimartingale,
while in $\rm (ii)$ we  only suppose the existence of $X$ (without regularity), but in
return we assume the existence of a smooth solution $G$ of the SPDE 
\eqref{eq: composition formulaA}.\\
 Otherwise, remark that the first line in  \eqref{eq: composition formulaA}, associated with $Y$, is purely the SDE part
of the dynamics of $G$ while the second corresponds to the partial differential part. This writing suggests a simple method
for solving such equations. Indeed, when considering any SPDE, if we are able to rewrite it in  the form   \eqref{eq:
composition
formulaA}, then we can hope solve it (if regularity of the identified coefficients holds) by associating two SDEs, and then
by  composing  with the solution associated with the partial differential part ($X$ in our result); this is the aim of
assertions $\rm{((iii),(iv)})$.  
\begin{proof} \rmi The proof is based on the It\^o-Ventzel formula, applied to $Y$
 as random field, and $\xi(t,z)$ as semimartingale; that leads to
\begin{eqnarray*}
dG(t,z)&=&\mu^Y(t,G(t,z))dt+\big(\sigma^Y(t,G(t,z))-Y_{y}(t,\xi(t ,
z))\xi_z(t, z)\sigma^X(t,z)\big).dW_t\\
 &+&\demi Y_{yy}(t,\xi(t,z))
\|-\xi_z(t,z)\sigma^X(t,z)\|^2+X_x(t,\xi(t , y))
{\widehat L}_{t,z}^Z(\xi) dt\\
&+& <dY_y(t,y),d\xi(t,z)>|_{y=\xi(t,z)}\\
&=&\mu^Y(t,\xi(t,z))dt+\big(\sigma^Y(.,G)-Y_{y}(.,\xi)
\xi_z\sigma^X \big)(t,z).dW_t \\
 &+&\big[\demi
Y_{yy}(.,\xi) \|\xi_z\sigma^X\|^2+Y_{y}(.,G)
{\widehat L}^X(\xi)-G_z\sigma^Y_y(.,G).\sigma^X
\big](t,z)dt
 \end{eqnarray*}
Now using identity $\partial_z\big(Y_y\big(t,\xi(t,z)\big)\big)=Y_{yy}(t,\xi(t,z))
\xi_z(t,z)$  and $
G_{z}(t,z)=Y_{y}(t,\xi(t,z))\xi_z(t,z)$, it 
follows, at first, that 
\begin{eqnarray*} 
 Y_{yy}(.,\xi) \|-\xi_z \,\sigma^X\|^2+ Y_{y}(.,\xi)
\partial_z\big(\|\sigma^X\|^2\xi_z\big)
&=& \partial_z\big(Y_{y}(.,\xi)\big)\big(\xi_z
\|\sigma^X\|^2\big)+Y_{y}(.,\xi)
\partial_z\big(\|\sigma^X\|^2\xi_z\big)
\big)\\
&=&\partial_z\big(Y_{y}(.,\xi)\xi_z \|\sigma^X
\|^2\big)=\partial_z\big(\|\sigma^X\|^2G_z\big)
\end{eqnarray*}
Second, by injecting this identity in the dynamics of $G$, we obtain
 \begin{eqnarray*}
 dG(t,z)&=&\big(\sigma^Y_t(G(t,z))-G_z(t,z)\sigma^X_t(
z)\big).dW_t\\ \nonumber
 &+&\Big(
{\widehat L}_{t,z}^X(G)+\mu^Y_t(G(t,z))-\partial_z\big(\sigma^Y(t,G(t,
z)).\sigma^X(t,z)\big) \Big)dt
\end{eqnarray*}
In a simpler formulation, that is equivalent to,
 \begin{eqnarray*}
dG(t,z)&=&\sigma^Y(t,G(t,z)).dW_t+\mu^Y(t,G(t,z))dt\\
 &-& \partial_zG(t,z)\sigma^X(t,
z)\big[dW_t+\sigma^Y_y(t,G(t,z)]+ {\widehat L}^{X}(G)(t,z)dt
\end{eqnarray*}
In the particular case, where $Y(t,x)=F(x)$, $\mu^Y\equiv 0$ and
$\sigma^Y\equiv 0$ and the result is obvious. \\
\rmii Again, It\^o-Ventzel  calculus yields to
\begin{eqnarray*}
d\xi(t,\overline{X}_t)&=&\Big[\frac{1}{2}\partial_y\big(||\sigma^X(t,y)||^2\xi(t,y)\big)(t,\overline{X}_t)-\mu^X(t,
\overline{X}_t)\xi_z(t,\overline{X}_t)\Big ]dt\\
&+&\xi_z(t,\overline{X}_t)\overline{\sigma}(t,\overline{X}_t)dW_t+\Big[\xi_z(t,\overline{X}_t)\mu^X(t,\overline{X}_t)+\frac{1
}{2}\xi_{zz}(t,\overline{X}_t)||\overline{\sigma}(t,
\overline{X}_t)||^2\\&-&\big(\xi_{
yy }
(t,\overline{X}_t)\sigma^X(t,\overline{X}_t)+\xi_{y}(t,\overline{X}_t)\sigma_z^X(t,\overline{X}_t)\big).\overline{\sigma}(t,
\overline{X}_t)\Big ] dt\\
&-&\xi_z(t,\overline{X}_t)\sigma^X(t,\overline{X}_t)dW_t
\end{eqnarray*}
and by arranging the terms properly we get
\begin{eqnarray*}
d\xi(t,\overline{X}_t)&=&-\xi_z(t,\overline{X}_t)\Big(\sigma^X(t,\overline{X}_t)-\overline{\sigma}(t,\overline{X}
_t)\Big)dW_t-\xi_z(t,\overline{X}_t)\big(\mu^X(t,\overline{X}_t)-\overline{\mu}(t,\overline{X}_t)\big)dt\\
&+&\frac{1}{2}\xi_{zz}(t,\overline{X}_t)\Big(||\overline{\sigma}(t,\overline{X}_t)||^2+||\sigma^X(t,\overline{X}
_t)||^2-2\sigma^X(t,\overline{X}_t).\overline{\sigma}(t,\overline{X}_t)\Big)dt\\
&+&\xi_z(t,\overline{X}_t)\Big(\sigma^X(t,\overline{X}_t)-\overline{\sigma}(t,\overline{X}_t)\Big).\sigma_z^X(t,\overline{X}
_t)dt\\
&=&-\xi_z(t,\overline{X}_t)\Big(\sigma^X(t,\overline{X}_t)-\overline{\sigma}(t,\overline{X}_t)\Big)dW_t-\xi_z(t,\overline{X}
_t)\big(\mu^X(t,\overline{X}_t)-\overline{\mu}(t,\overline{X}_t)\big)dt\\
&+&\frac{1}{2}\xi_{zz}(t,\overline{X}_t)||\overline{\sigma}(t,\overline{X}_t)-\sigma^X(t,\overline{X}_t)||^2dt\\
&+&\xi_z(t,\overline{X}_t)\Big(\sigma^X(t,\overline{X}_t)-\overline{\sigma}(t,
\overline{X}_t)\Big).\sigma_x^X(t,\overline{X}_t)dt
\end{eqnarray*}
Subsequently, by introducing new coefficients $ \Delta{\mu}$ and $  \Delta{\sigma} $ as follows:
$$ \Delta{\mu}(t,y):=\mu^X(t,y)-\overline{\mu}(t,\overline{X}_t) ~~\text{et}~~  \Delta{\sigma}(t,y):=
\sigma^X(t,y)-\overline{\sigma}(t,\overline{X}_t)$$
 the last line becomes 
$\partial_y\big(|| \Delta{\sigma}(t,y)||^2\xi_z\big)(t,\overline{X}_t)$
which yields to the desired formula:
\begin{eqnarray*}
d\xi(t,\overline{X}_t)&=&-\xi_z(t,\overline{X}_t)\Big(\sigma^X(t,\overline{X}_t)-\overline{\sigma}(t,\overline{X}
_t)\Big)dW_t+{\widehat
L}^{ \Delta}(\xi)(t,\overline{X}_t)dt\\
&=&-\xi_z(t,\overline{X}_t)\hat{\sigma}(t,\overline{X}_t)dW_t+{\widehat L}^{ \Delta}(\xi)(t,\overline{X}_t)dt
\end{eqnarray*}
\rmiii The properties of the random field $G(t,\overline{X}_t)$
 are obtained once again from It\^o-Ventzel's formula, in a very similar way than for $\xi(t,\overline{X}_t)$, since the
second line of
Equation \eqref{eq: composition formulaA} has the same form than the SPDE of $\xi$ except that the Brownian $dW$ is replaced
by $dW_.+\sigma^Y_x(G) dt$. So we obtain Equation \eqref{eq: compcharacteristicA}.
 \end{proof}
\paragraph{Back to  Progressive Utilities}
 Let now come back to progressive utilities ${\bf U}$.  The results
of this section will be of great use in the rest of this work,
especially when we focus on dynamic optimal portfolios. 
We give sufficient conditions on
progressive dynamic utilities so that assumptions of Theorem \ref{UxSDE} are
satisfied, in particular that the inverse $-\tU_y$ of $U_x$ is a semimartingale. 
These  assumptions are made on the coefficients of the intrinsic SDE$(\mu, \sigma) $  and not on the local characteristics
$(\beta, \gamma)$, since essential results are obtained from  SDEs properties.
\begin{theo}\label{Resume} Consider SDE$(\mu,\sigma)$ of  class $
\S^{1,\delta}, \delta\in (0,1]$ and let $Z$ be its unique monotonic solution of class  $\K^{1,\eps}_{loc}$ for any
$0<\eps<\delta$ .
For any deterministic utility function $u$  s.t. $u_x$ is
integrable near to $x=0$, define  $U(t,x)=\int_0^xZ_t(u_x(z))dz$.  Then\\[1mm]
\rmi $U(t,x)$ is an It\^o semimartingale with local characteristics $\beta(t,x)=\int_0^x\mu(t, Z_t(u_x(z))dz$ and
$\gamma(t,x)=\int_0^x\sigma(t,Z_t(u_x(z))dz$. \\[1mm]
\rmii  $U$ is a  progressive utility with derivative
$U_x(t,x)=Z_t(u_x(x))$. Moreover, $U$ is  a $\K^{2,\eps}_{loc}$-semimartingale for any $0<\eps<\delta$, with 
local characteristics $(\beta,\gamma)$ are of  class  ${\K}^{2,\eps}_{loc}\times\overline{{\K}}^{2,\eps}_{loc}$
for any $\eps<\delta$.\\[1mm]
\rmiii If  the SDE$(\mu,\sigma)$ is of  class $ \S^{m,\delta}, (m\ge 3, \delta\in (0,1])$, then, 
the progressive convex conjugate  utility $ \tU$ of $U$ is 
a $\K^{m-1,\eps}_{loc}$-semimartingale. \\
$-$ Its derivative $ \tU_y=-(U_x)^{-1}$ is a  $ \K^{m-2,\eps}_{loc}$-semimartingale, solution of the SDE$(\tilde
\sigma, \tilde \mu)$ in $\S^{m,\eps}$ for any $0<\eps<\delta$, where
\begin{equation}\label{eq:dualSDE}
\tilde \sigma(t,-z)=\displaystyle{\frac{\gamma_x(t,z)}{U_{xx}(t,z)}},\quad
\tilde \mu(t,-z)=
\frac{1}{U_{xx}(t,z)}
\Big(\beta_x(t,z)-\frac{1}{2}\partial_x\big(\frac{\|{\gamma}_x(t,z)\|^2}{\,U_{xx
}(t,z)}\big)\Big)
\end{equation}
$-$ Assumptions of Theorem \ref{ResPrA} are satisfied and the dynamics of
$\tU$ and $\tU_y$ are,
\begin{eqnarray*}
d\tU(t,y)&=&\gamma(t,-\tU_y(t,y)).dW_t+\beta(t,-\tU_y(t,y))dt +\frac{1}{2}
\tU_{yy}(t,y)\|\sigma_t(y)\big)\|^2\,dt \\
d\tU_y(t,y) &=&-\tU_{yy}(t,y)\sigma(t,y).dW_t+\Big(\demi\partial_y\big(\|\sigma\|^2\partial_y(\tU_y)\big)(t,
y)-\mu(t,y)\partial_y(\tU_y)(t,y))\Big)dt
\end{eqnarray*}
\end{theo}
 \noindent {\bf Concluding remarks on these two sections} Having introduced the progressive utilities and their convex
conjugate in  Section 1,
we have studied in detail  conditions ensuring  concavity and
 Inada conditions, by using that  $U_x$ is a monotonic solution of a "regular" SDE. 
But the main result given in  Section 2 concerns the inverse flow $\xi$ of the monotonic solution $X$ of a SDE, and the
sufficient conditions on
the coefficients  ensuring the semimartingale property of  the inverse flow $\xi$. The introduction of the adjoint SPDE
satisfied by $\xi$ and the result of Theorem \ref{pp: composition formula} are a powerful tool for that follows. We turn now
to the study of additional condition called   consistency property that we introduce
in the next section.

\section{Consistent Dynamic Utilities}\label{CDU}
\subsection{Definition of Consistent Dynamic Utilities}
The notion of progressive utility is very general and  should be  specified so
as to represent more realistically the dynamic evolution of the individual
preferences of an investor in a given financial market. The utility input
provides a differential constraint on the risk attitude of the investor in terms
of preference for higher or lower wealth, through the local risk tolerance
$\tau^U(x,t)=-\frac{U_x(x,t)}{U_{xx}(x,t)}$, and, that we call, 
 utility market risk premium $\eta^U(t,x)=\frac{\gamma_x(t,x)}{U_x(x,t)},$
revealing the interplay between the investment universe and the risk attitude. 

 The market input is described by a vector space $\GX$ of portfolios
incorporating feasibility and trading constraints and high liquidity. 
 Several interpretations of the subclass of admissible portfolios can
be done. The first one, proposed by Musiela and Zariphopoulou  is that this class 
describes all  investment universe. The second one is as follows: the market inputs may be viewed as a calibration
universe, and the class $\GX$ as a  test-class of processes. The existence of an
admissible portfolio giving the maximal satisfaction to the investor, which will
be preserved at all times in the future, explains  {\em the martingale property} in
the definition below.  On the other hand if the strategy in $\GX$ fails to be
optimal then it is better  not to  make investment. The optimal portfolio may be
viewed as a benchmark for the investor using the utility $U$.
 Once his consistent progressive utility  is defined,  an investor can then turn
to a portfolio optimization  problem  in a larger financial market or  to
calculate indifference prices. 
Following   \cite{zar-07a,zar-07}, a $\GX$-consistent dynamic utility  is defined as follows.
\begin{df}[$\GX$-consistent dynamic utility]\label{defUF}
A $\GX$-consistent dynamic utility ${\bf U}=\{U(t,x); t\geq 0, x>0\}$  is a
progressive utility
with the following additional properties:
\begin{description}
\item  {\sc Consistency with the test-class:}  For
any admissible wealth process  $X \in \GX$, 
\begin{equation*}
 \mathbb{E}(U(t,X_t) /{\cal F}_s)\leq
U(s,X_s), \>~ \forall s\leq t~~ 
a.s. 
\end{equation*}
\item{\sc Existence of optimal wealth:} 
For any initial wealth $x>0$, there exists an
 optimal wealth process $X^*\in \GX$ such that $X^*_0=x$, and for all $s\leq t$,
 \begin{equation*}
U(s,X^{*}_{s})=\mathbb{E}(U(t,X^{*}_{t})/{\cal F}_s) \>~ \forall s\leq t~~ 
a.s. 
\end{equation*}
\end{description}
 In short for any admissible wealth $X\in \GX$,
$U(.,X_.)$ is a positive supermartingale and a martingale  for the
optimal-benchmark wealth $X^*$.
\end{df}
\begin{Req}{\rm  \rmi  The martingale property  can be weakened by the following
localization procedure, if there exists a sequence of increasing stopping times
$T_n(X_0) \nearrow \infty$  on the random interval $[0,T_n(X_0)] $,
$U(.,X^*_.)$ is a martingale.\\
\rmii Note that the initial condition $U(0,.)=u(.)$ is part of the definition,
in particular $u$ is a deterministic utility function which  is fixed a priori
independently of  the  given financial market and in particular of $\GX$.\\
\rmiii So a $\GX$-consistent dynamic utility is a constant in time deterministic function only when  the
test-portfolios are local martingales. In this case, the optimal strategy is to
do nothing.\\
\rmiv Our definition differs slightly from the original (\cite{zar-07a})
 since  we do not  require that the wealth
processes $X$ are discounted.
This variation offers more options and allows us to study the invariance of the
class of dynamic utilities by change of numeraire. In any case, there is no
fixed horizon.}
\end{Req}

\subsection{The investment universe}
We consider a incomplete It\^o market, equipped with a $n$-standard Brownian motion,
$W$ with Brownian coordinates $(W_1,W_2,...,W_n)^T$  $(n\ge d)$ and characterized by a short rate $(r_t)$ and a
$n$-dimensional risk premium vector 
$(\eta)$. All these processes are defined on the filtered probability space $(\Omega,{\cal F}_{t\ge 0},\mathbb{P})$
satisfying usual assumptions, with minimal integrability assumptions. Since we only need to know the class of admissible
portfolios, we immediately give the mathematical definition of this class, based on the  self-financing equation without
arbitrage. The market incompleteness is modeled as in Lehoczky, Karatzas, Shreve \& Xu \cite{Karatzas04}.
The notations are the same than in
Karatzas and Shreve \cite{KaratzasShreve:01} where the interested reader may be find all complementary information.

\begin{df}[Test portfolios] \label{TestP} \rmi A positive It\^o semimartingale
$X^\kappa$ is called a portfolio, or admissible wealth process if
\begin{equation}\label{eq:DynamX}
dX^\kappa_t =X^\kappa_t  \big[r_t dt+\kappa_t. (dW_t + \eta_t dt)\big],~
\kappa_t \in \sigR_t.
\end{equation}
where $\kappa$ is a $n$-dimensional vector, progressive such that $\int_0^T\|\kappa_t\|^2dt<\infty, a.s.$, measuring the
volatility
vector of the wealth $X^\kappa$. \\
\rmii The family of admissible wealth processes, also called test portfolios is defined by some restrictions on the
volatility vector $\kappa$; we assume there exists progressive family of vector spaces $(\sigR_t)$  such that for any $t$, 
$\kappa_t\in \sigR_t$. \\
\rmiii The family of test portfolios is denoted by $\GX$. It may be easy to verify that $\GX$ is convex.
 \end{df}
\noindent The following short notations will be used extensively. 
Let $\sigR$ be a vector subspace of  $\R^n$. For any $x\in \R^n$,
$x^\sR$ is the orthogonal projection of the vector $x$
onto $\sigR$ and $x^{\perp}$ is the orthogonal projection
onto $\sigoR$. \\
The existence of  a risk premium  $\eta$ is a possible
formulation of the absence of arbitrage opportunity. Since from
\eqref{eq:DynamX}, the impact of the risk premium on the wealth dynamics only
appears through the term $\kappa_t. \eta_t $ for  $\kappa_t\in \sigR_t$, there is
a "minimal" risk premium $(\eta^\sR_t)$, the projection of $\eta_t$  on
the space $\sigR_t$ $(\kappa_t. \eta_t =\kappa_t. \eta^\sR_t )$,
to which we refer in the sequel. Moreover, the existence of $\eta^\sR$ is not
enough to insure the existence of equivalent martingale measure, since in
general we do not know if the exponential local martingale
$L^{\eta^\sR}_t=\exp(\int_0^t -\eta^\sR_s.dW_s-\frac{1}{2} \int_0^t
|\eta^\sR_s|^2\,ds)$ is a uniformly integrable martingale, density of an
equivalent martingale measure. Nevertheless, we are interested into the class of
the so-called state price density processes $Y^\nu$ (taking into account the
discount factor) who will play the same role for the dynamic conjugate utility,
than the wealth processes $X^\kappa$ for the dynamic  utility $U$.
\begin{df}[State price density process] \label{SPDP} \rmi A It\^o semimartingale
$Y^\nu$ is called a state price density process if  for any wealth process
$X^\kappa,~~\kappa \in \sigR$, $Y^\nu X^\kappa$ is a local martingale. It
follows that $Y^\nu$ satisfies,
\begin{equation}\label{Ynu}
dY^\nu_t=Y^\nu_t[-r_tdt+ (\nu_t-\eta^\sR_t).dW_t],\quad
\nu_t \in \sigoR_t, \quad Y^\nu_0=y
\end{equation}
\noindent
\rmii Denote $\GY$ the convex family of all state density processes $Y^\nu$ where $~\nu\in
\sigoR$ and observe that
 $Y^\nu$ is the product of $Y^0$ $(\nu=0)$ by the density martingale
$L^{\nu}_t=\exp\big(\int_0^t \nu_s.dW_s-1/2\int_0^t |\nu_s|^2ds\big)$.
\end{df}
\subsection{Consistent Dynamic Utility and Portfolio Optimization}\label{SPDECSPU} 
\paragraph{$\GX$-consistent dynamic utility and HJB constraint}
In Paragraph \ref{IPU}, more precisely in  Theorem \ref{UxSDE}, we have  characterized
progressive utilities in terms of their local characteristics 
$(\beta, \gamma)$ as well as in terms of the parameters $(\mu, \sigma)$ of the intrinsic SDE \eqref{eq:ItoSDE} satisfied by
$U_x$. In this section, we are concerned with the constraint induced on the drift characteristic $\beta$ of the dynamic
utility by the consistency property.
The consistency property  plays the same role that the dynamic programming
principle in the classical theory of backward expected utility maximization,
(see for example H. Pham \cite{Pham}). Thanks to It\^o-Ventzel's formula
(Theorem \ref{IVF}), constraints on the local characteristics $(\beta,
\gamma)$ of $\bf U$ lead to non standard Hamilton-Jacobi-Bellman Stochastic PDE. As in the classical case, the main
parameters of the SPDE are the risk tolerance process and the utility market risk premium.
\begin{df}[Utility risk tolerance and risk premium]\label{def:riskprem}
In this financial framework,  the utility risk tolerance random field is defined
by $\tau^U(t,x)=-\frac{U_{x}(t,x)}{U_{xx}(t,x)}$ and  the utility risk premium
random field by $\eta^U(t,x)=\frac{\gamma_x(t,x)}{U_x(t,x)}$ with its two
components $\eta^{U,\sR}\in \sigR, ~\eta^{U,\perp}\in \sigoR.$ \\[2mm]
{\rm Observe that Condition \eqref{eq: gamma condition} in Corollary  \ref{eq:Lipschitzcarac}  states that
$\eta^U$ is bounded in $x$ with random bound.}
\end{df}
\noindent The supermartingale property of $U(.,X^\kappa)$ implies that the drift of these processes must be
negative for all $\kappa \in  \sigR$, and equal to $0$ for some $\kappa^*$. 
We proceed by verification as in the classical case.
\begin{theo}[Utility-SPDE]\label{EDPGa}
Let $\bf U$ be a progressive utility which is  a $\Kc^{2,\delta}_{loc}$-semimartingale  $(\delta\in
(0,1])$ with local
characteristics $(\beta,\gamma)$. Assume the  drift constraint to be of
HJB type,    
\begin{eqnarray}\label{BETA}
\beta(t,x)=-U_x(t,x)r_t x -\frac{1}{2}U_{xx}(t,x)\inf_{\kappa \in
\sigR}\Big\{\|x \kappa\|^2+2 x \kappa.\big(\frac{U_x(t,x)\eta^\sR_t
+ \gamma_x(t,x)}{\,U_{xx}(t,x)}\big)\Big\}.
\end{eqnarray}
\rmi The   minimum of the quadratic form \eqref{BETA} is achieved at the optimal
policy $\kappa^*$  given by
\begin{equation}\label{xkappa}
\left \{
\begin{array}{cllll}
x\kappa^*_t(x)&= -\frac{1}{ U_{xx}(t,x)}(U_x(t,x)\eta^\sR_t
+\gamma_x^\sR(t,x))\\
\mbox{and }\>\> \beta(t,x)&=-
U_x(t,x)xr_t+\frac{1}{2}U_{xx}(t,x)\|x\kappa^*(t,x)\|^2
\end{array}
\right.
\end{equation}
\rmii For any $\kappa \in \sigR$, the process $U(.,X^{\kappa}_.)$ is a
supermartingale, and a local martingale for any solution (if there exists) $X^*$
of the SDE $dX^*_t=X^*_t\big(r_t dt +\kappa^*(t,X^*_t).(dW_t+\eta^\sigma_tdt)\big)$. Under additional integrability
assumptions, the $\GX$-consistency property is satisfied.
\end{theo}
\begin{Req} {\rm  In the classical backward framework,  similar SPDE is
investigated by Mania and Tevzadze\cite{Mania} using BSDE tools, and by Englezios and Karatzas \cite{KarEng}. (See also
Remark \ref{req:Karatzas}})
\end{Req}
\begin{proof}
 By It\^o-Ventzel's formula (Theorem \ref{IVF}), for any admissible portfolio
$X^\kappa$,
\begin{eqnarray*}
dU(t,X^\kappa_t)&=&\Big(U_x(t,X^\kappa_t)X^\kappa_t\>
\kappa_t+\gamma(t,X^\kappa_t)\Big).dW_t\\
&+&\Big(\beta(t,X^\kappa_t)+U_x(t,X^\kappa_t )r_t
X^\kappa_t+\frac{1}{2}U_{xx}(t,X^\kappa_t)\mathcal{Q}(t,X^\kappa_t,
\kappa_t)\Big)dt,\\
\mbox{where}\>\mathcal{Q}(t,x,\kappa)&:=&\|x \kappa\|^2+2 x
\kappa.\big(\frac{U_x(t,x)\eta^\sR_t +
\gamma_x(t,x)}{\,U_{xx}(t,x)}\big).
\end{eqnarray*}
Since $\kappa \in \sigR$, $\mathcal{Q}(t,x,\kappa)$ is only depending on 
$\gamma_x^\sR(t,x)$, the orthogonal projection of $\gamma_x(t,x)$ on
$\sigR_t$. 
The minimum $\mathcal{Q}^*(t,x)= \inf_{\kappa \in \sigR}\mathcal{Q}(t,x,\kappa)$
of the quadratic form $\mathcal{Q}(t,x,\kappa)$ is achieved at the optimal
policy $\kappa^*$  given by
\begin{equation}\label{kappaQ*}
\left\{
\begin{array}{cll}
x\kappa^*_t(x) &=-\frac{1}{\,U_{xx}(t,x)}\big(U_x(t,x)\eta^\sR_t
+\gamma_x^\sR(t,x)\big)\\[1mm]
\mathcal{Q}^*(t,x)&=-\frac{1}{U_{xx}(t,x)^2}\|U_x(t,
x)\eta^\sR_t+\gamma_x^\sR(t,x))\|^2=\,-\|x\kappa^*_t(x) \|^2.
\end{array}
\right.
\end{equation}
Then the drift of the semimartingale $U(t,X^\kappa_t)$ satisfies
\begin{eqnarray*}
&&\beta(t,X^\kappa_t)+U_x(t,X^\kappa_t )r_t
X^\kappa_t+\frac{1}{2}U_{xx}(t,X^\kappa_t)\mathcal{Q}(t,X^\kappa_t,\kappa_t) \\
&\leq& \beta(t,X^\kappa_t)+U_x(t,X^\kappa_t )r_t
X^\kappa_t+\frac{1}{2}U_{xx}(t,X^\kappa_t)\mathcal{Q}^*(t,X^\kappa_t,\kappa_t)\\
&=&\beta(t,X^\kappa_t)+U_x(t,X^\kappa_t )r_t
X^\kappa_t-\frac{1}{2}U_{xx}(t,X^\kappa_t)\|X^\kappa_t\kappa^*_t(X^\kappa_t)
\|^2
\end{eqnarray*}
The proof is complete.
\end{proof}

\paragraph{Conjugate of consistent dynamic utility }
The characteristics of the conjugate progressive utility $\bf \tU$ 
can be computed directly from Theorem \ref{ResPrA}. Given  that
$\beta$ is associated with an optimization program, we show that the dual drift
$\tilde{\beta}$ is also constrained by  a HJB type relation in the new
variables. So, the
convex conjugate utility $\bf \tU$ is  consistent with a family of state price
density processes (Definition \ref{SPDP}). As observed in Theorem \ref{Resume} (ii), the study of the conjugate utility $\tU$
requires stronger assumptions than the study of $U$.
\begin{theo}\label{thEDPDuale}
Let $\bfU$  a progressive utility with characteristics $(\beta, \gamma)$ satisfying Assumptions of Theorem \ref{Resume}. 
Then  its progressive convex conjugate utility $\bf \tU$ and its marginal conjugate utility ${\bf \tU_y}$ are It\^o random
fields with local characteristics $(\tilde\beta, \tilde\gamma)$ and $(\tilde\beta_y, \tilde\gamma_y)$ respectively.
Assume the  drift constraint of $\bfU$ to be of HJB type
\eqref{BETA}.
\\
 \rmi The local characteristics of the convex conjugate $\bf \tU$  are given by:
 \begin{eqnarray}\label{EDPSDuale'}
\left\{
\begin{aligned}
&~\tilde{\gamma}(t,y):=\gamma(t,-\tU_y(t,y)), \quad
\tilde{\gamma}_y(t,y):=-\gamma_x(t,-\tU_y(t,y)). \tU_{yy}(y)\\
&\tilde{\beta}(t,y)=y\tU_{y}(t,y)
r_t+\frac{1}{2\tU_{yy}(t,y)}\big(\|\tilde{\gamma}_y(t,y)\|^2-\|\tilde{\gamma}
^\sR _y(t,y)+y\tU_{yy}(t,y) \eta^\sR_t\|^2\big) 
\end{aligned}
\right .
\end{eqnarray}
\rmii The non linear drift $\tilde {\beta}(t,y)$ is associated with  the
following optimization program: 
\begin{eqnarray}\label{tildebeta}
\begin{aligned}
&\tilde {\beta}(t,y)=y\tU_y(t,y)r_t-\frac{1}{2}y^2\tU_{yy}(t,y)\inf_{\nu_t \in
\sigoR}\{
\|\nu_t-\eta^\sR_t\|^2+2\big(\nu_t-\eta^\sR_t\big).\big(\frac{\tilde{
\gamma}_y(t,y)}{y\tU_{yy}(t,y)} \big)\}
\end{aligned}
\end{eqnarray}
\rmiii The minimum of this quadratic form is achieved at the optimal policy 
 \begin{equation*}
y\nu^*_t(y)=
\frac{-\tilde{\gamma}^{\perp}_y(t,y)}{\tU_{yy}(t,y)}=\gamma^\perp_x\big(t,
-\tU_y(t,y)\big)=y\eta^{U,\perp}(t,-\tU_y(t,y)).
\end{equation*}
\rmiv Rewritten in terms of optimal strategy, $\tilde{\beta}(t,y)$ becomes
 \begin{equation}\label{eq:tildebetaoptimal}
\tilde{\beta}(t,y)=y\tU_{y}(t,y) r_t-\frac{1}{2}\tU_{yy}(t,y)\Big[\|y\nu^*_t(y)-
y\eta^\sR_t\|^2-2
(y\nu^*_t(y)- y\eta^\sR_t).\frac{\tilde{\gamma}_y(t,y)}{
\tU_{yy}(t,y)}\Big].
\end{equation}
\rmv For any admissible state price density process $Y^\nu\in \GY$ with $\nu\in \sigoR$, $\tU(t,Y^\nu_t)$ is a
 submartingale, and a local martingale for any solution $Y^*$ (if there exists) of the equation $dY^*_t=Y^*_t[-r_t
dt+(\nu^*(t,Y^*_t)-\eta^\sR_t).dW_t]$.
 \end{theo}
\begin{proof}\rmi
By Theorem
\ref{ResPrA}, the local characteristics $(\tilde{\beta},\tilde{\gamma})$ of
the conjugate random field $\tU$ are given by
$\tilde{\gamma}(t,y)=\gamma(t,-\tU(t,y))$ and
$\tilde{\beta}(t,y)=\beta^1(t,-\tU(t,y))$
where $\beta^1(t,x)=\beta(t,x)-\frac{1}{2U_{xx}(t,x)}\|\gamma_x(t,x)\|^2.$
Combining this identity with the HJB-constraint $\beta(t,x)=-
U_x(t,x)xr_t+\frac{1}{2}U_{xx}(t,x)\|x\kappa^*(t,x)\|^2$
 yields to
\begin{equation*}
\beta^1(t,x)=-xU_x(t,x)r_t-\frac{1}{2U_{xx}(t,x)}\Big(\|\gamma_x(t,
x)\|^2-\|U_x(t,x)\eta^\sR_t+\gamma_x^\sR(t,x)\|^2\Big).
\end{equation*}
(ii) \& (iii) Since the norm of the projection on $\sigR_t$ is the distance to
the orthogonal vector space $\sigoR_t$, 
$$\|\frac{\gamma_x^\sR(t,x)}{U_x(t,x)}+\eta^\sR_t\|^2:=\|\eta^{U,\sR}
(t,x)+\eta^\sR_t\|^2=
\inf_{\nu\in \sigoR_t}\|\nu-(\eta^{U,\sR}+\eta^\sR_t)\|^2.$$
Using the relation $|x|^2-|y|^2=|x-y|^2+2(x-y).y$, we get:
\begin{eqnarray*}
\|\eta^{U,\sR}(t,x)+\eta^\sR_t\|^2-\|\eta^{U}(t,x)\|^2&=&
\inf_{\nu \in \sigoR}
\{\|\nu-\eta^\sR_t\|^2+2(\nu-\eta^\sR_t).\eta^{U}(t,x)\}\\
&=&\|\eta^{U,\perp}(t,x)-\eta^\sR_t\|^2+2(\eta^{U,\perp}(t,x)-\eta^\sR
_t).\eta^{U}(t,x).
\end{eqnarray*}
By coming  back to $\beta^1$, we can make  the minimization program to be
explicit,
\begin{eqnarray*}
\beta^1(t,x)+xU_x(t,x)r_t&=&
\frac{U_x^2(t,x)}{U_{xx}t,x)}  \big(\inf_{\nu \in
\sigoR}\{\|\nu-\eta^\sR_t\|^2+
2(\nu-\eta^\sR_t).\eta^{U}(t,x)\}\big)\\
&=&\frac{U_x^2(t,x)}{U_{xx}t,x)}
\big(\|\eta^{U,\perp}(t,x)-\eta^\sR_t\|^2+2(\eta^{U,\perp}(t,
x)-\eta^\sR_t).\eta^ { U}(t,x)\big).
\end{eqnarray*}
The minimum in the quadratic form is achieved at $\eta^{U,\perp}(t,x)$,
corresponding  the optimal strategy $\nu^*_t(y)
=\eta^{U,\perp}(t,-\tU_y(t,y)).$ 
 From this and  the  identities,
$$\frac{U_x^2(t,-\tU_y(t,y))}{2U_{xx}(t,-\tU_y(t,y))}=-\frac{1}{2}y^2\tU_{yy}(t,
y),\quad \eta^{U}(t,-\tU_y(t,y))=-\frac{\tilde{\gamma}_y(t,y)}{y
\tU_{yy}(t,y)}$$
we get the desired formula for  $\tilde{\beta}$ both in  \eqref{EDPSDuale'} 
and in \eqref{eq:tildebetaoptimal}, i.e.
\begin{eqnarray*}
\tilde{\beta}(t,y)&=&y\tU_{y}(t,y)
r_t+\frac{1}{2\tU_{yy}(t,y)}\big(\|\tilde{\gamma}_y(t,y)\|^2-\|\tilde{\gamma}^\sR_y(t,y)+y\tU_{yy}(t,y) \eta^\sR_t\|^2\big)\\
&=&y\tU_{y}(t,y)
r_t-\frac{1}{2}\tU_{yy}(t,y)\Big[\|y\nu^*_t(y)-y\eta^\sR_t\|^2-2(y\nu^*_t(y)-y\eta^\sR_t).\frac{\tilde{\gamma}_y(t,y)}{
\tU_{yy}(t,y)}\Big].
\end{eqnarray*}
\rmiv is a simple rewriting of the HJB constraint on $\tilde \beta$.\\
\rmv The supermartingale property is proved in very similar manner than for dynamic utility, (proof of Proposition
\ref{thEDPDuale} (ii) and (iii)).
\end{proof}

\section{Marginal Utility SPDE and Optimal SDEs}\label{USPDEOSDE}
As seen in Section 1 and 2, our utility characterization is based on the marginal utility. So, we focus now on the
dynamics of the random field $\bf U_x$.
As established in Theorem \ref{EDPGa}, the drift characteristic of $\X$-consistent utility is constraint by a HJB type
condition \eqref{BETA} \\[1mm]
\centerline{$\beta(t,x)=- U_x(t,x)xr_t+\frac{1}{2}U_{xx}(t,x)\|x\kappa^*_t(x)\|^2.$}
The diffusion characteristic $\gamma_x$ is explained by the optimal policies of the primal and dual problems, Theorem
\ref{EDPGa} and
Theorem \ref{thEDPDuale}:  \\[1mm]
\centerline{$\gamma_x(t,x)=-U_{xx}(t,x)x\kappa^*_t(x)- U_x\eta^{\sR}_t+U_x(t,x)\nu^*_t(U_x(t,x))$}
 \\[1mm]
The characteristic $\gamma_x(t,x)$ may be rewritten in a more convenient form, using the diffusion coefficient of the optimal
 policy $\sigma^*_t(x):=x\kappa^*_t(x)$, and  the diffusion coefficient of the optimal state price density $ \tilde
\sigma^*_t(y)=- y \eta^{\sR}_t+\gamma^\perp_x(t,-\tU_y(t,y))$ as,  \\[2mm]
\centerline{$\gamma_x(t,x)=-\partial_x(U_{x})(t,x)\sigma^*_t(x)- U_x(t,x) \eta^{\sR}_t+  \tilde\sigma^*_t(U_{x}(t,x))$}
\\[1mm]
It is easy to recognize the diffusion coefficient of the SPDE  \eqref{eq: composition formulaA} in
Theorem  \ref{pp: composition formula} associated with the SDEs with diffusion parameters $\sigma^*_t(x)$ and $\tilde
\sigma^*_t(y)$. Moreover, by taking the $x$-derivative in the drift constraint $\beta(t,x)=-
U_x(t,x)xr_t+\frac{1}{2}U_{xx}(t,x)\|x\kappa^*_t(x)\|^2=- U_x(t,x)xr_t+\frac{1}{2}U_{xx}(t,x)\|\sigma^*_t(x)\|^2$,   it
appears 
naturally a divergence term associated with the optimal policy  $\sigma^*_t(x)=x\kappa^*_t(x)$ in the drift characteristic of
$\bf U_x$ which suggests to use the main theorem (Theorem \ref{pp: composition formula}) 
applied to the optimal SDEs.

\subsection{Main result : solving marginal utility SPDE via optimal SDEs}
To be closer to the notation of Theorem \ref{pp: composition formula}, we recall all the coefficients of SDEs  associated
with
the optimal policies $X^*$ and $Y^*$ {\em if they exist}:
\begin{equation}\label{eq:optcoeff}
\left\{
\begin{array}{llll}
\mu^*_t(x)&:=r_t x+x\kappa^*_t(x).\eta_t^\sR,\quad \sigma^*_t(x):=x\kappa^*_t(x) \\
\tilde \mu^*_t(y)&:=-r_t \,y,  \hspace{23mm} \tilde \sigma^*_t(y)=-\eta^\sR_t y+\gamma_x^\perp(t,-\tU_y(t,y))\\
\widehat L^*_{t,x}&:=\frac{1}{2}\partial_x(\|\sigma^*_t(x)\|^2 \partial_x)-\mu^*_t(x)\partial_x
\end{array}
\right.
\end{equation}
We start with the identification of the SPDE satisfied by the marginal  utility of a consistent dynamic utility imposing only
regularity condition on the utility random field and its local characteristics. We will then give additional conditions that
guarantee the existence of solution to  SDEs with coefficients $(\mu^*, \sigma^*)$ and $(\tilde \mu^*,\tilde \sigma^*)$. 
\begin{pp}\label{ppU_x}
 Let $U$ be a $\K^{2,\delta}_{loc}\cap\Cc^3$-regular $(\delta>0)$ progressive utility $U$,
whose the local characteristics $(\beta,\gamma)$ satisfy  the HJB constraints,
\begin{equation}\label{eq:HJBConst}
\left\{
\begin{array}{llll}
\gamma_x(t,x)&:=-U_{xx}(t,x)\sigma^*_t(x)+\eta^\sR _t U_x(t,x)+\tilde \sigma^*_t(U_x(t,x)) \\
\beta(t,x)&:=-U_x(t,x)x\,r_t+\frac{1}{2}U_{xx}(t,x)\|\sigma^*_t(x)\|^2
\end{array}
\right.
\end{equation}
  The marginal utility $\bf U_x$ is a decreasing solution of the SPDE\eqref{eq:
composition formulaA}  with coefficients $(\mu^*, \sigma^*)$ and $(\tilde \mu^*,\tilde \sigma^*)$\\[-7mm]
\begin{eqnarray}\label{eq: composition formulaU}
 dU_x(t,x)&=&\tilde \sigma^*_t(U_x(t,x)).dW_t+\tilde \mu^*_t(U_x(t,x))dt \nonumber\\
 &-&\partial_xU_x(t,x)\sigma^*_t(x).\big(dW_t+\tilde{\sigma}^*_y(t,U_x(t,x))dt)+\wL^*_{t,x}(U_x)dt
\end{eqnarray}
\end{pp}

\begin{proof}At first, as $U$ is assumed to be $ \K^{2,\delta}_{loc}\cap \Cc^{(3)}$-regular, $U_x$ is of class
$\K^{1,\delta}_{loc}$ and
its
local characteristics $(\beta_x,\gamma_x)$ are of class  $\Cc^1$ in $x$; then, the vectors $
\sigma^*_t(x)=x\kappa^*_t(x)=-(\gamma_x^\sR(t,x)+\eta^\sR _t U_x(t,x))/U_{xx}(t,x)$ and $\tilde \sigma^*_t(y)=-\eta^\sR_t
y+\gamma_x^\perp(t,-\tU_y(t,y))$ are also of class $\Cc^1$, necessary condition to define $\widehat L^*$.
By derivation of  the local characteristics, it is clear that $\beta_x$  contains a second order term in divergence
form associated with $x\kappa^*_t(x)=\sigma^*_t(x)$, which leads us  to
introduce the adjoint operator $\widehat L^*_{t,x}$ associated with $(\mu^*,\sigma^*)$.
 It remains to make some slight transformations on the drift characteristic. Observe that 
 \begin{equation}
\begin{array}{llll}
\beta_x(t,x)&=- \partial_x(U_x(t,x)xr_t)+\partial_x\big(\demi U_{xx}(t,x)\|\sigma^*_t(x)\|^2\big)\\
&=\widehat L^*_{t,x}(U_x)-r_tU_x(t,x)+\partial_xU_x(t,x) \sigma^*_t(x).\eta_t^\sR\\
&=\widehat L^*_{t,x}(U_x)+\tilde \mu^*_t(U_x) +\partial_xU_x(t,x) \sigma^*_t(x).\eta_t^\sR
\end{array}
\end{equation}
Let us give another interpretation of $ \sigma^*_t(x).\eta_t^\sR$. Since  $\tilde
\sigma^*(t,y)+\eta^\sR_t y$ belongs to the vector space $\sigoR_t$ 
 the spatial derivative $\tilde \sigma^*_y(t,y)+\eta^\sR_t$
 is also in $\sigoR_t$, yielding to the relation on the scalar products $-
\sigma^*_t(x).\eta_t^\sR=\sigma^*_t(x).\tilde\sigma^*_y(t,y)$. Then, Identity \eqref{eq: composition formulaU} holds true.
\end{proof}

 \noindent {\bf Comment} \noindent \rmi There is a fairly subtle relation between the SDE$(\mu^U,\sigma^U)$ introduced in
Theorem \ref{UxSDE}  to characterize the marginal utility $\bf U_x$  
 where $\mu^U(t,z)=\beta_x\big(t,-\tU_y(t,z)\big)$ and $\sigma^U(t,z)=\gamma^U_x\big(t,-\tU_y(t,z)\big)$, and the
coefficients of the two optimal SDEs, in particular in terms of diffusion coefficients. 
 The HJB constraint $\gamma^U_x(t,x)=-U_{xx}(t,x)\sigma^*_t(x)+\eta^\sR _t U_x(t,x))+\tilde \sigma^*_t(U_x(t,x)) $
becomes\\[-2mm]

$$\sigma^U(t,z)=\frac{ \sigma^*_t(-\tU_y(t,z))}{\tU_{yy}(t,z)}  +\eta^\sR _t z +\tilde \sigma^*_t(z) $$

 In particular, since $\tilde \sigma^*_t(z) =\sigma^{U,\perp}(t,z)$, any regularity property on $\sigma^{U}$ are immediately
transferred by linear projection on $\tilde \sigma^*_t(-z)$. But, only some local regularity on $\sigma^*_t(t,z)$ may be
deduced from global regularity of $\sigma^U$. Nevertheless, we can justify the existence of a global optimal wealth solution
in a similar way than for the inverse process. \\
\rmii Observe also that we obtained a way to generated $\X$-consistent utility only from their local characteristics and the
SPDE \ref{eq: composition formulaU},  since we do not use a priori the concavity assumption in the following theorem.

\begin{theo}[Main theorem]\label{maintheo}
Let $U$ be a  $\K^{2,\delta}_{loc}\cap \Cc^3$-semimartingale {\bf progressive utility $U$},
whose the local characteristics $(\beta,\gamma)$ satisfy  the HJB constraints \eqref{eq:HJBConst}; then, the derivative
$U_x(t,x)$
is solution of the SPDE \eqref{eq: composition formulaU}.\\
Assume, in addition, the existence of two positive adapted stochastic bounds $(K^1,K^2)$ such that
\begin{equation}\label{eq:gammaderivatives}
\|\gamma^\perp_x(t,x)\|\leq
K^1_t\,|U_x(t,x)|,  \> \|\gamma^\perp_{xx}(t,x)\|\leq K^2_t \,|U_{xx}(t,x)|, \>a.s., \> (K^1, K^2)\in \L^2(dt)
\end{equation}
{\sc Existence of optimal processes:} \rmi The  conjuguate SDE$(\tilde \mu^*, \tilde \sigma^*)$ is uniformly Lipschitz and
has a unique strong solution $Y^*_t(y)$, which is strictly positive, and strictly monotonic, with range $[0,\infty)$.\\
\rmii The SDE$(\mu^*, \sigma^*)$ has only local Lipschitz property, and admits a maximal  monotonic solution $X^*$ defined a
priori only up to a stopping times $\zeta^*(x)$. But,  the explosion time $\zeta^*(x)=\infty \, a.s. $ since the processes 
$U_x(.,X^*_.(x))$ is distinguishable from the solution $Y^*_.(u_x(x))$.\\
{\sc Consistency and marginal utility characterization:} The random field $U$ is a $\X$-consistent utility, with optimal
wealth $X^*$.  Furthermore, the derivative random field $U_x(t,x)$, solution of the SPDE \eqref{eq: composition formulaU}
is given by $U_x(t,x)=Y^*_t(u_x((X^*_t)^{-1}(x)))$. So $U_x(t,x)$ is a strictly decreasing and positive process with range
$[\infty,0]$.
\end{theo}
\begin{Req}{\rm  In this framework, we do not need that the inverse process $\X(t,x)$ is a semimartingale.
Nevertheless,  if $\gamma_x^\perp \equiv 0$, then $Y^*_t(y)=yY^0_t$ and as $U_x(t,x)$ is a semimartingale, $u_x(\X(t,x))$
is also semimartingale, consequently  $\X(.,x)$ is  a semimartingale since $U(t,x)$ and  $u$ are assumed to be of class
$\Cc^3$. }
\end{Req}
\begin{proof}
\rmi The assumption \eqref{eq:gammaderivatives} implies as above observed that the coefficients $(\tilde \mu^*, \tilde
\sigma^*)$ are uniformly Lipschitz in space, and then the SDE has an unique monotonic strong solution $Y^*_t(y)$. Since by
assumption $\|\gamma_x(t,x)\|\leq
K^1_t\,|U_x(t,x)|$, then $\| \tilde \alpha^*(t,y)\|=\frac{\|\tilde \sigma^*(t,y)\|}{y}\leq K^1_t $.
 As $K^1\in \L^2(dt)$, $Y^*$ which is also solution of the SDE $dY^*_t=Y^*_t[r_t dt+\tilde \alpha^*(t,Y^*_t).dW_t$ is
strictly positive.\\
\rmii Under the same assumption in addition to the property of $U$ of class $\K^{2,\delta}_{loc}$, the coefficients
$( \mu^*,\sigma^*)$ are only in ${\K}^{0,1}_{loc}\times \overline{{\K}}^{0,1}_{loc}\big)$, and then by Theorem \ref{SolMax},
there exists a  monotonic maximal solution
$X^*_t(x)$ up to explosion time $\zeta^*(x)$.\\
\rmiii a) We show together the property( $U_x(t,X^*_t(x))=Y^*_t(u_x))$ )and the fact that $\zeta^*(x)=\infty$. We are
applying, up to $\zeta^*(x)$, Theorem \ref{pp: composition formula} iv) to $U_x(t,x)$, solution of the SDPE \eqref{eq:
composition formulaU}, and to the process $X^*_t$ .Then, it follows that
 $U_x(t,X^*_t(x))$ is  solution of SDE$(\tilde \mu^*, \tilde \sigma^*)$ up to $\zeta^*(x)$ with initial condition $u_x$.
 In other words, $U_x(t,X^*_t(x))=Y^*_t(u_x(x))); t<\zeta^*(x)$.\\
 b)  When $t\mapsto  \zeta^*(x)$  on $\zeta^*(x)<\infty$, $U_x(t,X^*_t(x))\mapsto 0$, and $Y^*_t(u_x(x)))\mapsto
Y^*_{\zeta^*(x)}$, then $Y^*_{\zeta^*(x)}(u_x(x)))=0$, but  since $Y^*_t(y)$ is a strictly positive process, a contradiction
occurs if  $\P(\zeta^*(x)<\infty)>0$.\\
 \rmiv The property of $U$ are deduced of the identity $U_x(t,X^*_t(x))=Y^*_t(u_x(x)))$ and from Theorem \ref{Resume}.
\end{proof}
\begin{Req} \label{req:Karatzas}{\rm In \cite{KarEng} the authors have also shown that 
the solution of the backward SPDE can be represented through the composition  of two
invertible processes. The main difference with the approach proposed here is that
their processes  are  given as conditional expectation of monotonic functions  
as we describe them  in path-wise form as SDE solutions.  }
\end{Req}
\noindent{\bf Comment:} Obviously the main result, Theorem \ref{maintheo}, gives us only sufficient conditions to generate
$\X$-consistent utility  from their characteristics. In fact, the  assumptions gives us  the existence of a monotonic
strictly positive strong solution $Y^*$ of SDE$(\tilde{\mu}^*,\tilde{\sigma}^*)$, and by local regularity, the existence of a
monotonic solution $X^*$ of SDE$(\tilde{\mu}^*,\tilde{\sigma}^*)$ up to explosion time $\zeta^{x}$. Obviously, we can
interchange the roles of $Y^*$ and $X^*$. Indeed,  in 
 line with  Section 2, given the existence of a solution  of the SPDE \eqref{eq: composition formulaU}, by assertion $(iii)$
of Theorem \ref{pp: composition formula}, the existence of a solution  $X_.^*(x)$ starting from $x$ to the
SDE$(\mu^*, \sigma^*)$ implies that  $U_x(t, X^*_t(x))$ is solution of the SDE$(\tilde\mu^*, \tilde \sigma^*)$ and so is a
state price density process $\widehat Y^{*}_t(u_x(x))$, starting from $u_x(x)$.
Moreover, existence of a solution to SDE$(\mu^*, \sigma^*)$ is equivalent to existence of a solution
to SDE$(\tilde \mu^*,\tilde \sigma^*)$. Otherwise, if we give
ourselves $ X^* $  monotone solution of  SDE$(\mu^*,\sigma^*)$ with a semimartingale inverse $\X$ and $Y^*$ a 
$\K^2$ regular
solution to SDE$(\tilde{\mu}^*,\tilde{\sigma}^*)$ then, from assertion $(i)$ of Theorem \ref{pp: composition formula}, the 
compound process ${\widehat Y}^*_t\big(u_x\big(\X(t,x)\big)\big)$ is an obvious solution of SPDE \eqref{eq: composition
formulaU}.

\subsection{Reverse Engineering Problem}\label{RverseEnPr}
We consider the converse point of view of the marginal utility characterization given by Theorem
\ref{ppU_x}. More precisely, by taking as input an initial condition $U(0,.)=u(.)$ and some monotonic solution of 
SDE$(\mu,\sigma)$ and SDE$(\tilde{\mu},\tilde{\sigma})$,
 we propose an explicit way to recover all consistent  utility $U$ generating this wealth as optimal process. In the
classical expected
utility
framework, this reverse engineering problem has been considered by He and Huang  \cite{HeHuang} (1992) in a complete market.
Since the class of dynamic utilities is larger than the class of Markovian utilities considered in  \cite{HeHuang}, our
problem is easier to solve. In particular,  we establish that the only assumption  we need
is the monotonicity of the wealth process with respect to the initial wealth, plus some integrability conditions. \\
In the following theorem, we try to introduce  minimum assumptions used in  the different problems.
 Either we assume the  strong assumption that the SDE$(\mu,\sigma)$ is regular enough  to have a unique
monotonic strong solution $X$ whose the  inverse $\X$ is a semimartingale  or we make the existence of $X$  in addition to
the
existence of a solution $\X$ to the SPDE$(\wL^{\mu, \sigma}_{t,x}, -\sigma\partial_x)$.
\begin{theo}\label{ConvPv} Let $ \kappa\in \sigR$ be a volatility vector and  $\sigma_t(x)=x \kappa_t(x)$,
$\mu_t(x):=r_t x+\sigma_t(x).\eta_t^\sR$ the coefficients of the SDE defining an admissible wealth process. As
previously, let $\wL^{\mu, \sigma}_{t,x}:=\frac{1}{2}\partial_x(\|\sigma_t(x)\|^2 \partial_x)-\mu_t(x)\partial_x$ be its
adjoint
operator. \\
Similarly,  let $\nu\in \sigoR$ be  an orthogonal volatility vector, and $\big( \tilde{\mu}(t,y):=-yr_t,
\tilde{\sigma}(t,y)=y(\nu_t(y)-\eta^\sR_t)\big)$ the coefficients of the SDE defining an admissible state price density
process.\\
 {\sc Assumptions}
{\small
\rma {\sc Strong assumptions} Assume the SDE$(\mu,\sigma)$  in the class ${\cal S}^{3,\delta}$ $(\delta\in(0,1])$, 
so that SDE$(\mu,\sigma)$ has a unique monotonic solution $\bf X,$ whose the inverse $\bf  \X$ is solution of the SPDE
 \begin{equation}\label{SPDEIN}
d\X(t,x)=-\X_x(t,x)\sigma_t(x).dW_t+\widehat L_{t,x}(\X)dt.
\end{equation}

\rmb or {\sc Weak assumptions} Assume only the SDE$(\mu,\sigma)$  in the class ${\cal S}^{1,\delta}$ $(\delta\in(0,1])$,
and the existence of a solution $\bf \X$ of the SPDE$(\widehat L_{t,x}, -\sigma\partial_x)$

\rmc  Assume the SDE $(\tilde{\mu},\tilde{\sigma})$ in the class ${\cal S}^{2,\delta}$ with monotonic solution $Y$. }
 \\[1mm]
{\sc main Result}  \rmi For any  initial utility function $u$, 
the stochastic random field $\bf  V$ defined by $V(t,x)=Y_t\big(u_x\big(\X(t,x)\big)\big)$,  if it is integrable near
to zero, is the derivative of a
consistent stochastic utility $\bf U$, solution of the SPDE \eqref{eq: composition formulaU} with optimal wealth process 
${\bf X=\X^{-1}}$, solution of the SDE$(\mu,\sigma)$.\\
\rmii The derivative of the convex conjugate $\bf \tU$ of $\bf U$ is
${\bf \tU_y:=-V^{-1}}$ with $\tU_y(t,y)=X_t\big(-\tu_y\big(\Y(t,y)\big)\big)$ where $\bf \Y$ denote the inverse flow of $\bf
Y$.  Moreover, if the
SDE$({\tilde
\mu},{\tilde \sigma})$ belongs to $\S^{3,\delta}$,  the processes $(\tU_y(t,y))$ and $( \tU(t,y))$ are It\^o's
semimartingales. 

\end{theo}

\begin{proof}  Calculations are easy consequence of Theorem \ref{pp: composition formula} (ii),
applied to the processes $\X$ solution of the SPDE \eqref{SPDEIN}, and $Y$ in place of $X$ with initial condition $u_x$.
\\
\rmi With the notations of Theorem \ref{ConvPv}, we have 
\begin{eqnarray*}
 dV(t,x)&=&\tilde{\sigma}_t(V(t,x)).dW_t+{\tilde \mu}_t(V(t,x))dt\nonumber \\
 &-& \partial_xV(t,x)\sigma_t(x)\big[dW_t+{\tilde \sigma}_t(V(t,x))/V(t,x)dt+ {\widehat L}_{t,x}(V)dt
\end{eqnarray*}
Recall that ${\tilde \sigma}_t(y)=y(\nu_t(y)-\eta^\sR_t)$, so that ${\tilde
\sigma}_t(y).\sigma_t(x)=-y \eta^\sR_t.\sigma_t(x)$. So, the process $V$ satisfies the SPDE
\eqref{eq: composition formulaU}. As in the proof of Theorem \ref{ppU_x}, the drift $ \beta^V$ may be transformed into 
\begin{eqnarray*}
\beta^V(t,x)&=&{\tilde \mu}_t(V(t,x))-\eta^\sR_t.\sigma_t(x)V_x(t,x))+\partial_x\big(\demi
V_{x}(t,x)\|\sigma_t(x)\|^2\big)-\mu_t(x)V_x(t,x)\\
&=&- [r_t V(t,x)+V_x(t,x)(-x r_t+\eta^\sR_t.\sigma_t(x) -\eta^\sR_t.\sigma_t(x))+\partial_x\big(\demi
V_x(t,x)\|\sigma_t(x)\|^2\big)]\\
&=&- \partial_x(V(t,x)x\,r_t)+\partial_x\big(\demi V_x(t,x)\|\sigma(t,x)\|^2\big)
\end{eqnarray*}
Similarly, the diffusion characteristic of $V$, $\gamma^V$ is given by
\begin{eqnarray*}
\gamma^V(t,x)&=&\tilde{\sigma}_t(V(t,x))- V_x(t,x)\sigma_t(x)\\
&=&V(t,x)\nu_t(V(t,x))-\eta_t^{\sR}.V(t,x)-V_x(t,x)\sigma_t(t,x)
\end{eqnarray*}
\rmii We recognize that $V$ has the same local characteristics that marginal of consistent utility
(Theorem \ref{ppU_x} ). Taking the primitive of $V$, (if that makes sense) we define a random field $U(t,x)=\int_0^x
V(t,z)dz=\int_0^x Y\Big(t,u_x\big(\X(t,z)\big)\Big)dz$ which is a progressive utility satisfying the HJB constraint,
$\beta^U(t,x)=-U_x(t,x)xr_t+\frac{1}{2}U_{xx}(t,x)\|\sigma_t(t,x)\|^2$
and the diffusion constraint $
\gamma^U_x(t,x)=\gamma^V(t,x)=U_x(t,x)\nu_t(U_x(t,x))-\eta_t^{\sR}.U_x(t,x)-U_{xx}(t,x)x\kappa(t,x)$. So, $U$ is a
consistent dynamic utility.
\end{proof}
\begin{cor}
With the same notations as in Theorem \ref{ConvPv}, assume in addition global Lipschitz regularity on $(\mu,\sigma)$, that is
$\sigma\in \tilde{\K}^{0,1}_{\bf b} \cap \tilde{\K}^{3,\delta}_{loc} $,
$\nu\in\tilde{\K}^{0,1}_{\bf b} \cap \tilde{\K}^{2,\delta}_{loc} $ for $\delta\in (0,1]$. Let us also consider
a utility function $u$ satisfying Inada's conditions, such that $u_x\sim x^{-\zeta}$ $(\zeta<1)$ in the neighborhood of
$z=0$. \\
The composite random field $Y_t\big(u_x\big(\X(t,x)\big)\big)$ is
 integrable near to zero and it is the derivative of a
consistent stochastic utility $U$.
\end{cor}

\begin{proof}
Since $\sigma_t(x)=x \kappa_t(x)\in \tilde{\K}^{0,1}_b \cap \tilde{\K}^{3,\delta}_{loc} $,
$\nu_t(y)\in\tilde{\K}^{0,1}_b \cap \tilde{\K}^{2,\delta}_{loc} $ for $\delta\in (0,1]$  one can easily shows that the pair 
$(\mu, \sigma)$ belong to $\big(\K^{0,1}_b \cap
\K^{3,\delta}_{loc}\big)\times
\big(\tilde{\K}^{0,1}_b \cap \tilde{\K}^{3,\delta}_{loc}\big)$ and $(\tilde{\mu}, \tilde{\sigma})$  to $\big(\K^{0,1}_b
\cap \K^{2,\delta}_{loc}\big)\times \big(\tilde{\K}^{0,1}_b \cap \tilde{\K}^{2,\delta}_{loc}\big)$. Then, all assumptions of
Theorem \ref{ConvPv} are satisfied. Moreover, from Theorem \ref{theoSDE}, the unique strong monotonic solutions $X$ and
$Y$ satisfy the following asymptotic behavior (equation \eqref{limit2}) 
\begin{equation*}
 \lim_{z\rightarrow 0}\Big(\sup_{0\le t\le T}\frac{Z(t,z)}{z^{\eps}}\Big)=0  
\text{ and } 
\lim_{z\rightarrow 0}\Big(\sup_{0\le t\le
T}\frac{Z(t,z)}{z^{1+\eps}}\Big)=+\infty, \>\>\text{ for all }T
\end{equation*}
Consequently, one can easily shows  that, for any utility function $u$ satisfying Inada's conditions s.t. 
$u_x(x)<x^{-\zeta}$ for some $\zeta<1$, the composite random field $Y_t\big(u_x\big(\X(t,x)\big)\big)$ is also integrable
near to
zero. Indeed, it suffices to write that for any $\alpha,\beta,\gamma>0$ we have
\begin{eqnarray*}
\lim_{x\rightarrow 0}x^\alpha Y_t\big(u_x\big(\X(t,x)\big)\big)&=&\lim_{x\rightarrow 0}(X_t(x))^\alpha Y_t\big(u_x(x)\big)
=\lim_{y\rightarrow +\infty}(X_t(-\tu_y(y)))^\alpha Y_t(y)\\&=& 
\lim_{y\rightarrow +\infty}(\frac{X_t(-\tu_y(y))}{(-\tu_y(y))^\beta})^\alpha (-\tu_y(y)y^\gamma)^{\alpha \beta}
\frac{Y_t(y)}{y^{\alpha \beta \gamma}}
\end{eqnarray*}
As $u_x(x)\sim x^{-\zeta}$, we have 
 $\lim_{x\rightarrow 0}x^{\gamma'}u_x(x)=\lim_{y\rightarrow +\infty}-y^{\frac{1}{\gamma'}}\tu_y(y)=0$ for any
$\zeta<\gamma'<1$. Taking
$\gamma=\frac{1}{\gamma'}>1$, $\eps>0$ and $\beta=1+\eps$ we deduce, 
from the asymptotic behavior of $X^*$ and $Y^*$ and
Inada's conditions, that 
$$\lim_{x\rightarrow 0}x^\alpha Y_t\big(u_x\big(\X(t,x)\big)\big)
=\lim_{y\rightarrow +\infty}(\frac{X_t(-\tu_y(y))}{(-\tu_y(y))^\beta})^\alpha (-\tu_y(y)y^\gamma)^{\alpha \beta}
\frac{Y_t(y)}{y^{\alpha \beta \gamma}}=0,~\forall \alpha\in (0,\frac{1}{\beta \gamma})\subset (0,1),$$
which shows the integrability near to zero of $Y_t\big(u_x\big(\X(t,x)\big)\big)$.
 \end{proof}

\paragraph{Risk tolerance dynamics}
With the utility characterization given in Theorem \ref{ConvPv} the study of the risk
tolerance coefficient, taken along the optimal wealth, is greatly simplified. In
particular, the  martingale property established in He and Huang in
\cite{HeHuang} in a complete market, $Y^*=yY^0$ is easy to understand. A similar study has been performed by 
 Zariphopoulou \&  Zhou \cite{Zhou} in the special case of dynamic utility deduced by change of numeraire
(see next section) and stochastic environment from a deterministic time depending utility function. 
\begin{pp}\label{pprisk}
Assume  the optimal state density process to be  the minimal one $Y^*(.,y)=yY^0_.$, or equivalently the orthogonal risk
premium to be  $0$. Let $U$ be
a consistent dynamic utility with optimal regular wealth $X^*$ with derivative
$X^*_x.$\\
\rmi The path of the risk tolerance coefficient $\tau^U(t,x)=- \frac{U_x(t,x)}{U_{xx}(t,x)}$ at benchmark  optimal wealth is
given by $\tau^U(t,X^*_t(x)):=\tau^*(t,x)=\tau^u(x)X^*_x(t,x)$. 
\\
\rmii $\tau^U(t,X^*_t(x))=\tau^*(t,x)$ is an admissible portfolio with initial wealth $\tau^u(x)$ and admissible allocation
$\kappa^d_t=\kappa^*_t(X^*_t)+X^*_t\kappa^*_x(t,X^*_t)$. 
In particular, $Y^0_t \tau^*(t,x)$ is a local martingale.\\
\rmiii As in  \cite{Zhou}, the pair $(X^*, \tau^*)$  is solution of two dimensional SDE with random coefficients.
\end{pp}
\begin{proof} \rmi We start with the representation of the marginal utility as $U_x(t,x)=Y^0_t u_x\big(\X(t,x)\big)\big)$,
and of its derivative  $U_{xx}(t,x)=Y^0_t u_{xx}\big(\X(t,x)\big)\partial_x\X(t,x)\big)=Y^0_t
u_{xx}\big(\X(t,x)\big)/X^*_x(t,\X(t,x))$. By taking the ratio of these two quantities, we have that
 $$\tau^U(t,x)=\tau^u(\X(t,x))X^*_x(t,\X(t,x))\quad \text{or equivalently }\quad
\tau^U(t,X^*_t(x))=\tau^u(x)X^*_x(t,x)$$
\rmii This last characterization is interesting since from Equation \eqref{Z_zF}
\begin{equation}\label{eq:Xderive}
 dX^*_x(t,x)=X^*_x(t,x)\big[r_t dt +\sigma^*_x(t,X_t^*).(dW_t+\eta^\sR_t dt)\big]
\end{equation} 
As in the proof of Theorem \ref{ppU_x}, assertion $i)$, since $\sigR$ is a vector space, $X^*_x(t,x)$ is still an
admissible portfolio with allocation policy $\kappa^d_t=\sigma^*_x(t,X_t^*)=\kappa^*_t(X^*_t)+X^*_t\kappa^*_x(t,X^*_t)$. \\
\rmiii The pairs of processes $(X^*(t,x),X^*_x(t,x))$, and $(X^*(t,x),\tau^*(t,x))$ are solutions of the same 2-dimensional
SDE.
\end{proof}
\noindent All  results of the proposition can be extended to the general case, where $Y^*_t(y)$ is no more a linear function
of $y$,
but the interpretation of the results is more difficult.

\begin{pp} Let $U$ be a $\K^{2,\delta}_{loc}\cap\Cc^3 $-consistent dynamic utility $(\delta>0)$, then\\
\rmi The risk tolerance coefficient $\tau^U$ is given by
\begin{equation}\label{eq:risktolerance}
\tau^U(t,x)=\frac{Y^*_t(u_x(\X))X^*_x(t,\X)}{u_{xx}(\X)Y^*_x(t,u_x(\X))}=
\Big(\tau^u(\X)\frac{Y^*_t(u_x(\X)X^*_x(t,\X)}{u_{x}(\X)Y^*_y(t,u_x(\X))}
\Big)(t,x).
\end{equation}
and, along the optimal wealth,   $\tau^U(t,X^*_t(x))=\frac{Y^*_t(u_x(x))}{Y^*_y(t,u_x(x)) u_{xx}(x)}
X^*_x(t,x)$.\\
\rmii The derivative of the optimal wealth is an admissible portfolio associated with  the allocation
$\kappa^d_t=\kappa^*_t(X^*_t)+X^*_t\kappa^*_x(t,X^*_t)$ and initial wealth $1$.\\
\rmiii The process $Y^*_y(t,y)\tau^U(t,X^*_t(x))$ is a local martingale.
\end{pp}
The proof of this proposition is obvious.
%



\paragraph{Pathwise dynamic programming principle} \label{par:pdpp}In this paragraph we are interested
in using the  previous marginal utility characterization from a backward point of view. Unlike the standard case, we obtain a
pathwise dynamic programming principle. \\
For this, observe that the identity $ U_x(t,x)=Y^*_t\big(u_x(\X(t,x))\big)$ leads to the following identity based on  the
inverse flow $\Y$ of $Y^*$,
$ u_x(x)=\Y_t(U_x(t,X^*_t(x))\big).$
This explains how to recover the  marginal utility at time $0$ from the stochastic terminal one at a time $ t $. This point
of view is
interesting,  because in this case the initial condition $ u_x$ is simply the derivative of the value function of a
classical optimization problem. To go further into this idea, properties of stochastic flows, previously not used (as
explained above), will be essential.\\
{\bf Flow property}.

\rmi Let $ (Z _ t (x)) $ be a  strictly monotonic  random field with respect to $x$ with range $ [0, \infty) $ whose the
inverse random field is denoted by
 $ \Z_t(z)$. We extend the random field $\bf Z$ to intermediate dates $ (0\le s\le
t) $ by defining $Z_t(s,x)=Z_t(\Z_s(x))$. The following classical notation $Z_t(s,x):=Z_{s,t}(x)$ is useful to express the
semi-group property, that is for
$s\leq t\leq u$, $Z_{t,t}(x)=x$ and  $Z_{t,u}\circ Z_{s,t}=Z_{s,u}$. Note that $Z_t(x)=Z_{0,t}(x)$.\\
Similarly, we can extended the random field $\bf \Z$ to intermediate dates by taking the inverse of  $(Z_t(s,x))$ such that
$\Z_t(s,z)=Z_s(\Z_t(z)):=\Z_{s,t}(z)$. Then $\Z_t(s,z)=\Z_{s,t}(z)$ is a backward flow achieving the amount $z$ at date $s $.
Then, the semi-group property  holds true with the inverse order of the dates, that is
for $s\leq t\leq u$, $\Z_{s,t}\circ \Z_{t,u}= Z_s\circ \Z_t \circ Z_t\circ \Z_u=\Z_{s,u}$.\\

\rmii When the flow $Z_t(x)$ is the monotonic solution of  some SDE$(\mu,\sigma)$ starting from $x$ at time $0$, a classical
result states that $Z_{s,t}(x)$ is a solution defined on $[s,\infty)$ of the same SDE$(\mu,\sigma)$ starting from $x$ at time
$s$. When  the SDE$(\mu,\sigma)$ is regular enough so that $\Z_t(z)$ is the monotonic solution of the SPDE$(\mu,\sigma)$,
then $\Z_{s,t}(z)$ is  solution on $[s,\infty)$ of the same SPDE with initial condition $z$ at time $s$. Given that 
$\Z_{s,t}(z)=Z_s(\Z_t(z))$, this result may be viewed as a consequence of Theorem \ref{pp: composition formula} (i), applied
to the regular function $G(y)= Z_s(y)$ considered as deterministic after the time $s$. But, it is more natural to consider
$\Z_{s,t}(z)$ as a process in $s$
and so to introduce, time reversal and backward integration as in Kunita \cite{Kunita:01} Section 4.5, or Carmona \& Nualart
\cite{Carmona}. We do not develop this point of view here.

 Let us come back to consistent utility framework.
\begin{pp} 
We adopt the same framework as in Theorem \label{ConvPv} with the monotonic solutions $X$ and $Y$ of two SDE$(\mu,\sigma)$
and SDE$(\tilde \mu,\tilde \sigma)$. As usual the inverse processes are denoted by $\X$ and $\Y$.\\
The marginal utility is defined as $U_x(t,x)=Y_t\big(u_x(\X_t(x))\big)$.\\
\rmi For any $t$, the pathwise identity holds true, $u_x(x)=\Y_t\big(U_x(t,X_t(x))\big)$. This property is close to the 
following one, also true on the classical backward point of view,
\begin{equation} \label{eq:pathwiseid}
u_x(x)=\Y_t\big(U_x(t,X_t(x))\big), \quad  u_x(x)=\E\big(U_x(t,X(t,x))X_x(t,x)\big)
\end{equation} 
\rmii More generally,  using the stochastic flows $X^x_{s,t}, Y^y_{s,t}, \X^z_{s,t}, \Y^u_{s,t}$, we have the pathwise
dynamic principle
\begin{equation} \label{eq:pathwise PDP}
 U_x(t,x)=Y_{s,t}\big(U_x(s,\X_{s,t}(x))\big), \quad U_x(s,x)=\Y_{s,t}\big(U_x(t,X_{s,t}(x))\big),
\end{equation} 
\end{pp}

\begin{proof}
The both identities are easily deduced from the identity
$$u_x(\X_t(x))=\Y_t(U_x(t,x))=u_x(\X_{s}(\X_{t,s}(x))=\Y_s(U_x(s,\X_{t,s}(x))).$$
\end{proof}
\noindent {\bf Comment} The link between HJB SPDE and the inverse flows of some SDEs  has received little attention in the
literature to our knowledge, may be essentially because people are more interested by the Markovian case.

\section{Openness to other topics and works}\label{EDSU}
We close the paper by some openness to other topics and works;  we show  the stability of
the notion  of consistent utility by change of numeraire and then,  without loss of generality,  we can  consider the 
martingale market where the
portfolios are simple local martingales and the stochastic PDE's are easier to deal with.  
We also apply our method to the specific
example of decreasing consistent utilities (see \cite{Mike}  and \cite
{zar-08a}) where the volatility vector $\gamma$ is given equal to zero, given a new interpretation of the optimal wealth as
solution of inf-convolution problems in random power utilities.
\subsection{Change of numeraire}
 Notations about the  investment universe are the same as in Section \ref{CDU}.
A numeraire is a monetary reference used
as payment instrument. When  transactions take place in a given country,  the domestic currency is used as numeraire, 
but in an international setting, a common numeraire is used in general for all the transactions. Because an investor is
often
faced with investing in different markets, we are concerned in this paragraph, by the impact of the change of
 numeraire on its progressive utility random field.\\
It is well-known (see Geman, El Karoui, and Rochet \cite{Geman}) that the self-financing property is
invariant by change of numeraire. So, if $N$ is an It\^o positive continuous semimartingale, solution of linear equation
with
coefficients $({ \mu}^N x,{ \delta}^N x)$, it is well-known (see also Platen and Heath, or Karatzas and Kardaras,
\cite{Geman, Platen, Karatzas01}) that using $N$ as new numeraire
transforms an   It\^o market with  risk premium $\eta^\sR$ and short rate $r_t$ into an
investment universe, where the admissible portfolios are the processes ${\widehat X}^{\hat \kappa}= X^{\kappa}/N$ and the
state price density processes are the processes $\widehat Y^{\nu}=N.Y^{\nu}$.  The parameters of this new market are $\hat
\eta=\eta^\sR-\delta^N$, for the risk premium and  $\hat r=r- \mu^N + \delta^N.\eta^\sR$ for the short rate.\\
\rmi By the previous results, any ${\cal X}$-consistent dynamic utility $U$ defines a ${\cal X}^N$-consistent dynamic utility
$U^N$ by the transformation $U^N(t,\hat{x})=U(t,\hat{x}N_t)$. Then the conjugate $\tU^N$ is given by $\tU^N(\hat
y)=\tU(t,\hat{y}/N_t)$.\\
\rmii The class ${\cal X}^N$ of the admissible portfolios is characterized by the processes $\hat \kappa=\kappa-\delta^N$. 
The vector spaces $(\cal R_t; t\geq 0)$
are transformed into affine spaces $\widehat{\cal R}_t=\cal R_t-\delta^N_t$. Nevertheless if $\delta^N_t\in {\cal R}_t$, for
any time $t$, the constraint spaces $\widehat{\cal R}_t$ in the new market are the same than  the constraint spaces ${\cal
R}_t$ in the initial market.\\
\rmiii The associated optimal  portfolio is
$X^{N,*}=X^*/N$ and the optimal state price density process is $Y^{N,*}=Y^*/N$.\\
\rmiv  \,The diffusion characteristic  $ \gamma^N(t,\hat x)$ is
obtained from It\^o-Ventzel's formula, 
$ \gamma^N(t,\hat{x})=\gamma(t,\hat{x}.N_t)+\hat{x}U^N_{\hat{x}}(t,\hat{x})\delta^N_t$. The new drift $\beta^N(t,{\hat x})$
is more complicated to described, and the explicit form is not  of great interest in the general case.\\
\rmv In the usual case, the market numeraire $(Y^0)^{-1}=M$ (also called numeraire portfolio, or growth optimal
portfolio as in  Platen and Heath \cite{Platen}) is chosen as numeraire $Y^0$. Since $M$ is an admissible portfolio with
volatility
$\eta^{\sR}$,  the market parameters $\hat \eta$ and $\hat r$ are identically zero,
and  any admissible portfolio and any admissible state price density is a local martingale (under the historical
probability).\\
The optimal policies are simply given by 
$${\hat x}{\hat \kappa}^{*}_t({\hat x})=\frac{\gamma^{M,\sR}_{\hat x}(t,{\hat x})}{U^{M}_{{\hat x}{\hat x}}(t,{\hat
x})},\quad {\hat y}{\hat\nu^{*}}_t({\hat y})=\frac{{\tilde
\gamma}_{\hat y}^{M,\perp}(t,{\hat y})}{\tU^{M}_{{\hat y}{\hat y}}(t,{\hat y})}$$
while, the drift characteristics of the utility $U^M(t,{\hat x})$ and  $\tU^M(t,{\hat y})$ are given by,
 $$
\beta^M(t,{\hat x})=-\frac{1}{2}U^M_{{\hat x}{\hat x}}(t,{\hat x})\|{\hat x}\kappa^{*,M}_t({\hat x})\|^2,~
{\tilde \beta}^M(t,{\hat y})=-\frac{1}{2}{\tilde U}^M_{{\hat y}{\hat y}}(t,{\hat y})\|{\hat y}\nu^{*,N}_t({\hat y})\|^2,\quad
$$
Consequently, $U^M$ is a supermartingale and $\tU^M$ is a submartingale. Moreover, if $\gamma^M_{\hat x}\in \sigR$ the
conjugate utility 
$\tU^M$ is a local martingale and the optimal dual process is constant, $Y^{M,*}\equiv1$. By symmetry, if $\gamma^M_{\hat
x}\in\sigoR$, $U^M(.,\hat x)$ is a local martingale for any $\hat x$, and  the optimal wealth $X_t^{M,*}(\hat x)\equiv \hat
x.$

\subsection{Decreasing Consistent Utilities}  In this section, all prices of the investment universe are assumed to be
discounted, corresponding to the case where $r\equiv 0$.
An interesting class of consistent utilities is the class of decreasing consistent utilities,
which was studied  and fully characterized in the literature by Musiela \& al. \cite{zar-08a} and Berrier \& al. \cite{Mike}.
This utilities have a volatility characteristic $\gamma$ identically zero. It is an example where the  dual SPDE is easier to
study than the  primal one, since by taking $\gamma=0$, 
it follows from Theorem \ref{thEDPDuale} \eqref{EDPSDuale'}, that $U$ and $\tU$  are solutions of the following SPDEs
\begin{equation}\label{SPDEdecr}
dU(t,x)=\frac{1}{2}\frac{U_x(t,x)^2}{U_{xx}(t,x)}||\eta^{\sR}_t||^2 dt, \quad
d\tU(t,y)=-\frac{1}{2}y^2\tU_{yy}(t,y)||\eta^{\sR}_t||^2dt.
\end{equation}
which implies by convexity, that $t\mapsto U(t,.),~ \tU(t,y)$ are decreasing functions. \\[-6mm]
\paragraph {Example of power utilities}\rmi It is easy to verify that the power dual utility functions
$\tU^\vartheta(t,y)=\frac{1}{1-\vartheta}(1-{\tilde C}_t ^\vartheta y^{1-\vartheta}), \> (\vartheta>0)$, where as  usual the
parameter $\vartheta$ is the risk tolerance coefficient, are \label{'solutions of the dual PDE'} \ref{SPDEdecr},
if and only if  ${\tilde C}_t^\vartheta(\omega)$ is solution of the ordinary equation $d{\tilde C}^\vartheta_t (\omega)=- 
{\tilde C}_t ^\vartheta(\omega)  \vartheta(1-\vartheta) ||\eta^\sigma_t(\omega)||^2dt$.  In other words, 
$$\tU^\vartheta(t,y)=\frac{1}{1-\vartheta}(1-{\tilde C}_t ^\vartheta y^{1-\vartheta}),\> \, {\tilde C}^\vartheta_t
=\exp(-{\tilde \epsilon}(\vartheta)A^\eta_t),\>\mbox{with} \,A^\eta_t=\int_0^t||\eta^{\sR}_s||^2ds, \>\,{\tilde
\epsilon}(\vartheta)=(1-\vartheta)\vartheta>0.$$ 
Then, $\tU^\vartheta(t,yY^0_t)-\frac{1}{1-\vartheta}= -\frac{y^{1-\vartheta}}{1-\vartheta} (Y^0_t)^{1-\vartheta}{\tilde
C}^\vartheta_t $ is a  martingale, since $y Y^0_t$ is the optimal state price density. \\
\rmii Let us observe that  at any time $t>0$ the marginal conjugate utility $\tU^\vartheta_y(t,y)=-{\tilde C}_t ^\vartheta
y^{-\vartheta}$ is no longer a monotonic function of the risk tolerance coefficient $\vartheta$ since the function $\vartheta
\mapsto {\tilde \epsilon}(\vartheta)A^\eta_t+\vartheta \ln(y)$ is no monotonic.\\
\rmiii The marginal power utility $U^\vartheta_x(t,x)$ is given by $U^\vartheta_x(t,x)=({\tilde C}_t^\vartheta)^{1/\vartheta}
\,x^{-1/\vartheta}$. Since $(Y^0_t)^{1-\vartheta}{\tilde C}^\vartheta_t $ is a local martingale, the process 
$X^{*,\vartheta}_t(x)=x{\tilde C}^\vartheta_t (Y^0_t)^{-\vartheta}$ is an admissible portfolio, and it is easy to see that it
is the optimal one. Using the notation ${\overline X}^{\vartheta}_t=X^{*,\vartheta}_t(1)= {\tilde C}^\vartheta_t
(Y^0_t)^{-\vartheta}$, we have $X^{*,\vartheta}_t(x)=x{\overline X}^{\vartheta}_t$.\\
\rmiv As in the deterministic case, the optimal strategy $\kappa^{*,\vartheta}(t,x)$ is collinear to the risk premium
$\eta^{\sR}_t$ with factor $\vartheta$, $\kappa^{*,\vartheta}(t,x)=\vartheta\eta^\sR_t$. 
 Then, we recover the characterization of the marginal utility at the optimal state price density $yY^0_t$ as
$\tU^\vartheta_y(t,yY^0_t)=-{\tilde C}_t ^\vartheta
y^{-\vartheta}(Y^0_t)^{-\vartheta}=\tu^{\vartheta}_y {\overline X}^{\vartheta}_t.$ \\[-8mm]
\paragraph {Characterization of decreasing conjugate utilities} \rmi The set of positive solutions $\tU$ to the dual linear
PDE \ref{SPDEdecr} is a convex cone, stable by positive linear  combination. From this, it is natural to consider the
integral of conjugate power utilities, with respect to some positive  Borel measure $m$, including the initial condition. The
avoid the problem related to the constants, we formulate the problem on the marginal conjugate utilities and assume that $
\int_{\R^*_+}y^{-\vartheta}dm(\vartheta)<\infty$, to define the new conjugate marginal utility
\begin{equation}\label{dualrepres}
\tU_y^m(t,y)=\int_{\R^*_+}\tU^\vartheta_y(t,y){\tilde C}^\vartheta_t dm(\vartheta), \quad 
\tu_y^m(y)=-\int_{\R^*_+}y^{-\vartheta}dm(\vartheta)
\end{equation}
Such assumption on the initial condition $\tu_y(y)$ (then denoted $\tu^m_y(y)$) is equivalent to say that $\tu^m_y(e^z)$ is a
completely monotonic function.
All conjugate utility functions $\tu(y)$
do not verify this  condition, but if it the case and 
if $ m $ is compactly supported, it is easy  to check that any primitive $ \tU^m (t, y) $ verifies the random PDE
\eqref{SPDEdecr}, since the property is true for the power conjuguate functions $\tU^\vartheta_y(t,y)$ using Dirac measure at
$\vartheta$ as shown in page\,\pageref{'solutions of the dual PDE'}.
 In other words, $ \tU^m (t, y)(\omega) $
 is a space-time harmonic function of a geometrical Brownian  motion with  variance $ A_t(\omega) $. By a generalization of
Widder's Theorem \cite{Widder}, Musiela \& al. \cite{zar-08a} and  Berrier \& al. \cite{Mike} showed that there is
no other solution to the random PDE \eqref{SPDEdecr}.\\[-8mm]
 \paragraph{Sup-convolution interpretation} \rmi There is an interesting interpretation of these stochastic utilities: one
can imagine an investor starting with power
utility, with some ambiguity on its risk tolerance coefficient. At time $0$, she starts with a mixture of marginal conjugate
power utilities, weighted by some measure $m$. At time $t$, the marginal conjugate consistent utility is still a mixture of
power conjugate utilities with respect to the  measure $m_t(d\vartheta):= {\tilde C}^\vartheta_t\,dm(\vartheta)$, with
decreasing in time random density  ${\tilde C}^\vartheta_t.$
The stochastic measure $m_t(d\vartheta)$ is the unique measure which ensure that the process $\tU^m_y$ constructed in
equation \eqref{dualrepres} is the derivative of the conjugate of a consistent utility.\\
\rmii {\sc Optimal wealth} The characterization of the marginal conjugate utility $\tU_y^{m}(t,yY^0_t)$ along the optimal
density process $yY^0_t$
as $\tU_y^{m}(t,yY^0_t)=X_t^{*,m}(-\tu^m_y(y))$  is useful to characterize the optimal wealth $X_t^{*,m}(x)$, since the same property is true for the standard power utility functions,
\begin{eqnarray*}
X_t^{*,m}(-\tu^m_y(y))=\int_{\R^*_+}\tU^\vartheta_y(t,yY^0_t)dm_t(\vartheta)=
\int_{\R^*_+}X_t^{*,\vartheta}(-\tu^{ \vartheta}_y(y))dm_t(\vartheta)
= \int_{\R^*_+}\overline{X}_t^{\vartheta}y^{-\vartheta}dm_t(\vartheta),
\end{eqnarray*}
To give a decomposition directly in terms of $x$, we start  with a family of well-chosen initial wealths $x^{\vartheta,m}$
satisfying
$ x^{\vartheta,m}(x)=- \tu_y^{\vartheta}(u^m_x(x))$ so that  $\int_{\R^*_+}
x^{\vartheta,m}(x))dm(\vartheta)=-\tu^m_y(u^m_x(x))=x.$
So, the optimal wealth process $X^{*,m}_t$ issued from $x=-\tu^m_y(u^m_x(x))$ is given by the closed formula 
\begin{eqnarray}\label{CNRX*}
X^{*,m}_t(x)=\int_0^\infty X_t^{*,\vartheta}( x^\vartheta(x)) m_t(d\vartheta)=\int_0^\infty x^\vartheta(x)\overline{X}_t^{\vartheta}m_t(d\vartheta), \mbox{with} \int_0^\infty  x^\vartheta(x)m(d\vartheta)=x.
\end{eqnarray}
 $X^{*,m}_t(x)$  is strictly increasing and regular with respect to its initial condition $x$, since any function 
$x^{\vartheta,m}(x)$ is monotonic and regular.  The conjuguate of $\tU^{m}_t(y)$, $U^m(t,x)$ is a consistent utility with optimal wealth $X^{*,m}_t(x)$.\\
\rmiii \, {\sc The risk tolerance coefficient}  We are concerned with the risk tolerance coefficient $\displaystyle \tau_U^m(t,x)=-\frac{U^m_x(x)}{U^m_{xx}(x)}$ of the utility function $U^{m}(t, x)$, given that for power utility $\tau^{\vartheta}(t, x)=\vartheta x$.\\
According
to assertion $(i)$ of Proposition \ref{pprisk} stating that $\tau_U^m(t,X^{*,m}_t(x))=\tau_u^m(x)X^{*,m}_x(t,x)$, we see that 
\begin{eqnarray*}
\tau_U^m(t,X^{*,m}_t(x))=\tau_u^m(x)X^{*,m}_x(t,x)=\int_0^\infty \tau_u^m(x)x_x^\vartheta(x)X_x^{*,\vartheta}(t, x^\vartheta(x))
m_t(d\vartheta)
\end{eqnarray*}

On the other hand, by the fact that $x^\vartheta_x(x)= -u^m_{xx}(x)\tu^\vartheta_{yy}(u^m_x(x))$,  we easily get the equality 
$\tau^\vartheta(x^\vartheta(x))=x_x^\vartheta(x)\tau^U(x)$, and so
and hence,
\begin{eqnarray}
\tau^m_U(t,X^{*,m}_t(x))=\int_0^\infty\tau^\vartheta(t,X^{*,\vartheta}_t( x^\vartheta(x))) m_t(d\vartheta)
\end{eqnarray}
We still have some mixture properties along the optimal processes.\\

\rmiv \, {\sc Sup-convolution interpretation} Classical result in convex analysis shows the link between mixture of conjuguate utility functions and inf-convolution problem, in the following form: assume that the derivative of some conjuguate utility function $\tu^{m}_y (y)$ may be represented as the integral of some family of marginal conjuguate utility functions $\tu^{\vartheta}_y(y)$, with respect to some positive Borel measure $m$, 
$\tu^{m}_y (y)=\int_0^\infty  \tu^{\vartheta}_y(y) m(d\vartheta)$, 
where the integral is finite for any $y$. Then,  $\tu^{m}$ is the convex conjuguate of the sup-convolution problem associated with 
\begin{eqnarray} \label{supconv}
u^{m}(x)= \sup\Big\{ \int_0^\infty u^{\vartheta}(x_{\vartheta})m(d\vartheta)\big| (x_{\vartheta}) \mbox{ such that}\>\int_0^\infty  x_{\vartheta}m(d\vartheta)=x\Big\}
\end{eqnarray}
The optimal solution $(x^*_{\vartheta})$ that achieves the maximum in the optimization problem \ref{supconv} is given
explicitly as  
$x^{\vartheta,m}(x)=\tu_y^{\vartheta}(-u^m_x(x))$.\\
Come back to our problem of decreasing consistent utility random field. The strategy associated with the initial condition
 $x^{\vartheta,m}(x)$ may be interpreted as the deterministic optimal allocations of the initial wealth with respect to the parameter $\vartheta$, in a Pareto optimal equilibrium where a continuum of agents with different risk aversion are in competition. 
The same interpretation holds true at time $t$, with the family of wealth $X_t^{*,\vartheta}( x^\vartheta(x))$  and the random measure  $m_t(d\vartheta)$.\\
Moreover, the decreasing consistent utility random field is given at any time $t$ as,
\begin{eqnarray} \label{exaxtform}
U^{m}(t, x)= \sup\Big\{ \int_0^\infty U^{\vartheta}(t, X_{\vartheta}(t,x))m_t(d\vartheta)\big| (X_{\vartheta}(t,x) ) \mbox{ such that}\>\int_0^\infty  X_{\vartheta}(t,x)m_t(d\vartheta)=x\Big\}
\end{eqnarray}
\noindent
\rmv {\sc The role of the initial utility function} In the study of decreasing utilities, we introduced an assumption on the initial condition, that is
$\tu_y^m(y)=-\int_{\R^*_+}y^{-\vartheta}dm(\vartheta)$. The Borel measure $m$ is determining in the definition of the optimal wealth. This may seem at odds with the rest of the paper. This is not the case and even this is a nice example to illustrate our results.

Indeed, starting from the optimal portfolio
$X^{*,m}$ defining in Equation \eqref{CNRX*}, it is easy to construct a new consistent utility (more easily its Fenchel conjuguate ${\tilde V}(t,y)$) having $(X^{*,m}, yY^0)$  as optimal processes, starting from an initial concave function $v$ with Fenchel conjuguate $\tilde v$. Thanks to the dual characterization $\tV_y^{m}(t,yY^0_t)=X_t^{*,m}(-\tv_y(y))$ of the marginal conjuguate utility,  the same identity applied to $\tU_y^{m}$ yields to the following relation $\tV_y^{m}(t,yY^0_t)=X_t^{*,m}(-\tv_y(y))=\tU_y^{m}\big(t,u^m_x(-\tv_y(y))Y^0_t\big)$. It remains to make the change of variable $y\mapsto y/Y^0_t$, to obtain the marginal dual conjuguate
$$\tV_y^{m}(t,y)=X_t^{*,m}(-\tv_y(y/Y^0_t))=\tU_y^{m}\big(t,u^m_x(-\tv_y(y/Y^0_t))Y^0_t\big).$$
Obviously, by this transformation we lose the decreasing property in time of $\tV_y^{m}(t,y)$ since $u^m_x(-\tv_y(y/Y^0_t))Y^0_t$ is no more a decreasing process. The same kind of construction may be made when the optimal dual process $yY^0_t$ is replaced by a monotonic one.

\noindent
Note, similar ideas are developed in \cite{MRADNEK03} to build richer classes of utilities.\\[-8mm]

\paragraph{Conclusion}
In this new approach, the solution of the utility SPDE have a pathwise representation,
unlike to the characteristics method  where  the solutions are represented as  a conditional expectation. 
There are several advantages of this connection between  SPDEs and SDEs due to the many results of the SDE theory. 
To the best of our knowledge, there are no or few results that assert the monotonicity or the convexity of such solutions. Also, there may
be other advantages in numerical methods and simulations of the SDE than of SPDE.

Otherwise, this paper investigates consistent stochastic utilities from the
 SPDE point of  view. This leads therefore to make  strong
regularity assumptions: the market is a Brownian market and securities are
modeled as continuous semimartingales. Utilities are at least of class
$\Kc^{2}$ in the sense of Kunita in  order to apply It\^o-Ventzel's formula
and to deduce the SPDEs. Moreover, the method of stochastic utilities
construction   is based on the  dynamics of stochastic flows and their
inverses, and therefore additional regularity assumptions   on $X^*$ and
$Y^*$ are required. However, one can take a direct approach still based on  monotonicity assumptions on optimal processes
for the primal and dual problem, and on compound flows formula ;
it is showed in \cite{MRADNEK02}, that  these assumptions can
be considerably weakened. Indeed,  considering any financial market in
which the securities are modeled as  bounded semimartingales, the stochastic utilities  are of class $\Kc^{1}$ and wealth 
process are required to lie in a convex class $\GX \subset \PX$, the
 monotonicity assumption  of $X^*$ and $Y^*$ is sufficient to show the validity
of the construction  proposed in this work,  using analysis methods and
optimality conditions.\\

\nocite{*}
\bibliographystyle{plain}
\bibliography{ArefNewNov012}

\end{document}